\documentclass[a4paper,preprint,11pt]{elsarticle}

\usepackage[T1]{fontenc}
\usepackage{lmodern}
\usepackage{microtype}
\usepackage{amsmath,amssymb,amsthm,mathtools,bm}
\usepackage{enumitem}
\usepackage{booktabs}
\usepackage{array}
\usepackage{hyperref}
\usepackage[nameinlink,noabbrev]{cleveref}
\hypersetup{
  colorlinks=true,
  linkcolor=blue,
  citecolor=blue,
  urlcolor=blue,
  pdftitle={Dependence Across Birth Cohorts in Finite-Age Branching Processes},
  pdfauthor={Mojtaba Moradi}
}

\newtheorem{theorem}{Theorem}[section]
\newtheorem{proposition}[theorem]{Proposition}
\newtheorem{corollary}[theorem]{Corollary}
\newtheorem{lemma}[theorem]{Lemma}
\theoremstyle{definition}
\newtheorem{assumption}[theorem]{Assumption}
\newtheorem{definition}[theorem]{Definition}
\theoremstyle{remark}
\newtheorem{remark}[theorem]{Remark}

\newcommand{\E}{\mathbb{E}}
\newcommand{\Var}{\operatorname{Var}}
\newcommand{\Cov}{\operatorname{Cov}}
\newcommand{\cum}{\operatorname{cum}}
\newcommand{\1}{\boldsymbol{1}}

\newcommand{\Nzero}{\mathbb{N}_0}
\newcommand{\cF}{\mathcal{F}}
\newcommand{\fall}[2]{(#1)_{\underline{#2}}}

\journal{Stochastic Processes and their Applications}

\begin{document}
\begin{frontmatter}
\title{Joint Cohort Laws and Mixed Cumulants with Spectral Asymptotics for Growing Gaps in Finite-Age Branching Processes}

\author[aff1]{Mojtaba Moradi}
\address[aff1]{Department of Industrial Engineering, Faculty of Technology and Engineering-East of Guilan, University of Guilan, Rudsar, Iran\\
mmoradi@guilan.ac.ir}

\begin{abstract}
We investigate the joint distribution and temporal dependence of birth cohorts in a discrete-time age-dependent branching model in which an individual may reproduce at several integer ages. Under an independent age-specific reproduction array, we derive a recursive joint probability generating function, an exact history-dependent transition kernel, and a martingale innovation representation. These identities yield explicit covariance formulae and finite coefficient recursions for mixed factorial moments and mixed cumulants. When reproduction is restricted to a finite maximum age $L$, the cohort vector admits a finite-dimensional companion-matrix representation. In the primitive supercritical case, we establish an exact $L^p$ criterion: for every integer $p\geq2$, the Perron martingale and the normalized cohort state converge in $L^p$ if and only if all age-specific offspring distributions have finite $p$th moments. The normalized cohorts consequently share a common limit, which satisfies an age-indexed smoothing-transform equation with an explicit triangular cumulant recursion. Our principal asymptotic result concerns two cohorts whose separation increases with time. Under a simple-real-subdominant-root condition, we obtain a three-regime asymptotic equivalent for the mean-square distance between their mean-normalized sizes. The rate is governed by whether the squared modulus of the dominant non-Perron root is smaller than, equal to, or larger than the Perron root. This result quantifies the approach to the common cohort limit and distinguishes the present cohort-level question from first-order Perron convergence alone.
\end{abstract}

\end{frontmatter}

\section{Introduction}\label{sec:intro}

Branching processes describe populations whose members generate new individuals, but the timing of reproduction matters whenever individuals can reproduce more than once. In a one-generation Galton--Watson process, each individual contributes to the next generation once. In age-dependent models, by contrast, individuals born in different cohorts may contribute offspring to the same calendar-time cohort. The resulting sequence of birth-cohort sizes is therefore generally dependent and need not itself be a first-order Markov chain.

The broad framework of general branching processes, including Crump--Mode--Jagers processes, accommodates repeated reproduction at different ages and overlapping generations; a foundational general age-dependent formulation is due to Crump and Mode \cite{CrumpMode1968}. Discrete-time formulations and their random-tree representation are discussed by Jagers and Sagitov \cite{JagersSagitov2008}. The analysis of moments in age-dependent branching processes has a substantial history. Sevast'yanov and Chistyakov \cite{SevastyanovChistyakov1971} developed multidimensional renewal equations for branching-process moments. Mode \cite{Mode1974} derived discrete-time renewal equations for mean and covariance functions in a generalized one-type age-dependent model. Kimmel \cite{Kimmel1983} used a point-process formulation to obtain symbolic recursions for arbitrary finite-dimensional distributions in age- and time-dependent branching processes. More recently, Janson \cite{Janson2019} studied second-order fluctuations in lattice Crump--Mode--Jagers processes and identified distinct regimes governed by the roots of the associated renewal equation. Cheng \cite{Cheng2024} established limit and central limit results for supercritical age-structured branching models in which reproduction depends on an individual's remaining lifetime. That remaining-lifetime framework differs from the discrete repeated-age cohort model considered here, but it is relevant to the current development of age-structured branching limit theory.

The model considered here is a deliberately specific member of this broader family. An individual born at time $k$ may produce offspring at ages $1,2,\ldots$, and the distribution of the number produced at age $t$ depends on $t$ but not on calendar time or cohort size. We impose mutual independence of all age-specific reproduction variables, including variables corresponding to different ages of the same individual. A general Crump--Mode--Jagers reproduction point process need not satisfy this last condition. The present assumptions therefore yield a factorization of the conditional generating function by parental age, which is the structural feature used throughout the paper.

A closely related discrete-time period-dependent branching model was studied by Fathi Vajargah and Moradi~\cite{FathiMoradi2012}, following their earlier age-/period-dependent epidemic formulation~\cite{FathiMoradi2011}. In particular, the 2012 model already allows an individual to contribute offspring at more than one age and is therefore a direct predecessor of the population mechanism considered here. The present paper instead makes the independent reproduction array explicit and develops joint-cohort distributional identities, arbitrary-order mixed-moment and mixed-cumulant recursions, and finite-age asymptotic results that were not the focus of those earlier papers. The precise relation and remaining overlap are discussed in Section~\ref{sec:discussion}.

The dependence between birth cohorts is a central feature of this model. Although offspring contributions scheduled at a fixed calendar time are conditionally independent across parental ages, the cohort sizes that determine the numbers of summands are themselves dependent. This distinction rules out a general factorization of the joint cohort pgf into marginal pgfs. A simple calculation illustrates the point. Under the model assumptions,
\[
Z_1=X_{0,1}(1),\qquad
Z_2=\sum_{j=1}^{Z_1}X_{1,j}(1)+X_{0,1}(2).
\]
Since the two age-specific reproduction variables of the initial ancestor are independent and the reproduction variables of the first-born cohort are independent of $Z_1$, we have
\[
\E[Z_2\mid Z_1]=m_1Z_1+m_2,
\qquad
\Cov(Z_1,Z_2)=m_1\Var(Z_1)=m_1\sigma_1^2.
\]
Thus $Z_1$ and $Z_2$ are generally dependent. More generally, the total-variance decomposition contains cross-covariance terms:
\[
\begin{aligned}
\Var(Z_n)
={}&\sum_{t=1}^n\sigma_t^2a_{n-t}
+\sum_{t=1}^nm_t^2\Var(Z_{n-t})\\
&+2\sum_{1\leq t<u\leq n}m_tm_u\Cov(Z_{n-t},Z_{n-u}).
\end{aligned}
\]
These terms are essential whenever the earlier cohorts are dependent.

The paper has two complementary parts. The first develops exact finite-time distributional identities: a joint cohort pgf, a history-dependent transition kernel, an orthogonal innovation representation, and coefficient recursions for mixed factorial moments and cumulants. The second concerns finite maximum reproductive age. We identify the associated finite-dimensional companion recursion, establish a necessary-and-sufficient age-by-age moment condition for $L^p$ convergence of the Perron martingale, and characterize the common cohort limit through a smoothing transform. The main additional result is a quantitative theorem for growing cohort gaps. It gives a sharp three-regime mean-square asymptotic under an explicit simple-root expansion of the renewal sequence. The threshold is the square-root scale of the Perron root, and the constants are expressed in terms of the offspring variances and the full mean renewal sequence. We do not claim that the Perron martingale or smoothing-transform mechanism is new; the contribution is the explicit cohort-level formulation and, in particular, the growing-gap mean-square asymptotic.

The paper is organized as follows. Section~\ref{sec:model} specifies the reproduction array. Sections~\ref{sec:joint} and~\ref{sec:second} develop the joint law and second-order structure. Section~\ref{sec:allmom} gives arbitrary-order moment and cumulant recursions. Section~\ref{sec:asymptotic} contains the finite-age Perron theory, the common-limit cumulants, and the growing-gap theorem. Section~\ref{sec:discussion} discusses the relation to earlier work and the scope of the results.

\section{Model and conditional structure}\label{sec:model}

\subsection{Birth cohorts and reproduction variables}
Let $Z_n$ denote the number of individuals born at integer time $n$, with one initial ancestor:
\[
Z_0=1,\qquad Z_n\in\Nzero,\quad n\geq 1.
\]
For an individual $j$ born at time $k$, let $X_{k,j}(t)$ be the number of offspring that individual produces at age $t$. Potential individuals in each cohort are indexed by $j\geq1$; only indices $j\leq Z_k$ correspond to individuals actually present in cohort $k$.

\begin{assumption}[Independent age-specific reproduction]\label{ass:ind}
For each $t\geq1$, the random variables $\{X_{k,j}(t): k\geq0,\ j\geq1\}$ have a common distribution, denoted by the distribution of $X(t)$. The full family
\[
\{X_{k,j}(t): k\geq0,\ j\geq1,\ t\geq1\}
\]
is mutually independent, and each $X(t)$ takes values in $\Nzero$.
\end{assumption}

Thus age may affect the offspring distribution, but the law does not depend on calendar time, cohort size, or the reproduction of other individuals. In particular, the counts produced by the same parent at two different ages are independent. This is stronger than independence between distinct parents and is used explicitly in the product formulae below.

The birth-cohort recursion is
\begin{equation}\label{eq:population}
Z_n=\sum_{t=1}^{n}\sum_{j=1}^{Z_{n-t}}X_{n-t,j}(t),\qquad n\geq1.
\end{equation}
We use the conventions that an empty sum is zero and $Z_k=0$ for $k<0$. Define
\[
\phi_t(s)=\E[s^{X(t)}],\qquad m_t=\E[X(t)],\qquad
\sigma_t^2=\Var(X(t)),
\]
whenever the relevant moments exist.

\subsection{Filtration and the conditional generating function}
For a precise filtration, regard the variables $X_{k,j}(t)$ as preassigned for every potential index $(k,j,t)$, including indices $j$ that may not correspond to an individual ultimately born. Define
\[
\cF_n=\sigma\{X_{k,j}(t): k\geq0,\ j\geq1,\ t\geq1,\ k+t\leq n\}.
\]
Then $Z_0,\ldots,Z_n$ are $\cF_n$-measurable. By Assumption~\ref{ass:ind}, all variables scheduled at calendar time $n$ are independent of $\cF_{n-1}$. The use of potential variables for unborn individuals only completes the underlying independent array; the population recursion uses only indices $j\leq Z_k$.

\begin{proposition}[Conditional probability generating function]\label{prop:conditional-pgf}
Under Assumption~\ref{ass:ind}, for every $n\geq1$ and $s\in[0,1]$,
\begin{equation}\label{eq:conditional-pgf}
\E[s^{Z_n}\mid\cF_{n-1}]
=\prod_{t=1}^{n}\phi_t(s)^{Z_{n-t}}.
\end{equation}
\end{proposition}
\begin{proof}
Conditionally on $\cF_{n-1}$, the cohort sizes $Z_{n-t}$ are fixed. For each $t$, the contribution
\[
Y_{n,t}:=\sum_{j=1}^{Z_{n-t}}X_{n-t,j}(t)
\]
is a sum of $Z_{n-t}$ independent copies of $X(t)$, and therefore has conditional pgf $\phi_t(s)^{Z_{n-t}}$. For distinct $t$, the variables involved in these sums have distinct indices $(n-t,j,t)$ and are mutually independent by Assumption~\ref{ass:ind}. Their conditional pgfs multiply, proving \eqref{eq:conditional-pgf}.
\end{proof}

\begin{remark}
The conditional factorization in \eqref{eq:conditional-pgf} depends on independence across ages for a fixed parent. If the vector $(X_{k,j}(1),X_{k,j}(2),\ldots)$ had dependent coordinates, the age-wise factors could not in general be separated in this way. A formulation with a general reproduction point process would then be required.
\end{remark}

\section{Means and the joint cohort law}\label{sec:joint}

\subsection{Mean renewal equation}
\begin{theorem}[Mean recursion]\label{thm:mean}
If $m_t<\infty$ for $1\leq t\leq n$, then, with $a_n:=\E[Z_n]$,
\begin{equation}\label{eq:mean-rec}
a_n=\sum_{t=1}^{n}m_ta_{n-t},\qquad a_0=1.
\end{equation}
\end{theorem}
\begin{proof}
Taking conditional expectation in \eqref{eq:population}, or differentiating \eqref{eq:conditional-pgf} at $s=1$, gives
\[
\E[Z_n\mid\cF_{n-1}]=\sum_{t=1}^{n}m_tZ_{n-t}.
\]
Taking expectations yields \eqref{eq:mean-rec}.
\end{proof}

\begin{corollary}[Composition representation]\label{cor:composition}
For $n\geq1$,
\begin{equation}\label{eq:composition}
a_n=\sum_{k=1}^{n}\ \sum_{\substack{t_1+\cdots+t_k=n\\t_i\geq1}}
 m_{t_1}\cdots m_{t_k}.
\end{equation}
\end{corollary}
\begin{proof}
Repeated substitution in \eqref{eq:mean-rec} selects an ordered sequence of positive lags. The substitution terminates when the remaining time index is zero, at which point the terminal factor is $a_0=1$. The ordered sequences selected in this way are exactly the compositions of $n$, and each contributes the product in \eqref{eq:composition}.
\end{proof}

\begin{corollary}[Constant age-specific mean]\label{cor:constant-mean}
If $m_t=m$ for every $t\geq1$, then
\[
a_0=1,\qquad a_1=m,\qquad a_n=m(m+1)^{n-1},\quad n\geq1.
\]
\end{corollary}
\begin{proof}
For $n\geq1$, \eqref{eq:mean-rec} gives $a_n=m\sum_{j=0}^{n-1}a_j$. For $n\geq2$, subtracting the corresponding identity at $n-1$ yields $a_n=(m+1)a_{n-1}$. Since $a_1=m$, the stated formula follows. In particular, the initialization $a_0=1$ must be kept separate; $a_n=(m+1)^n$ is not the correct sequence for this model.
\end{proof}

\subsection{Joint probability generating function}
For $n\geq1$, define
\[
H_n(s_1,\ldots,s_n)=\E\left[\prod_{r=1}^{n}s_r^{Z_r}\right],
\qquad H_0=1.
\]
\begin{theorem}[Joint pgf recursion]\label{thm:joint-pgf}
For $n\geq1$,
\begin{equation}\label{eq:joint-pgf}
H_n(s_1,\ldots,s_n)=\phi_n(s_n)
H_{n-1}\big(s_1\phi_{n-1}(s_n),\ldots,
 s_{n-1}\phi_1(s_n)\big).
\end{equation}
\end{theorem}
\begin{proof}
By Proposition~\ref{prop:conditional-pgf},
\[
\begin{aligned}
H_n(\boldsymbol s)
&=\E\left[\left(\prod_{r=1}^{n-1}s_r^{Z_r}\right)
\E[s_n^{Z_n}\mid\cF_{n-1}]\right]\\
&=\E\left[\prod_{r=1}^{n-1}
\{s_r\phi_{n-r}(s_n)\}^{Z_r}\right]\phi_n(s_n),
\end{aligned}
\]
where the last factor corresponds to the single initial ancestor $Z_0=1$. The resulting expression is precisely \eqref{eq:joint-pgf}.
\end{proof}

\begin{remark}
The joint recursion retains the dependence among birth cohorts. In general, the marginal pgfs of $Z_1,\ldots,Z_n$ cannot be multiplied to obtain $H_n$, because the cohort sizes share ancestors and earlier reproduction events.
\end{remark}

\subsection{Finite-dimensional transition probabilities}
Let $p_t(x)=\Pr\{X(t)=x\}$, and let $p_t^{*q}$ denote the $q$-fold convolution of $p_t$, with $p_t^{*0}=\delta_0$.

\begin{proposition}[History-dependent transition kernel]\label{prop:kernel}
For $z_0=1$ and any history $(z_1,\ldots,z_{n-1})$ having positive probability,
\begin{equation}\label{eq:kernel}
\begin{aligned}
&\Pr(Z_n=z_n\mid Z_0=z_0,\ldots,Z_{n-1}=z_{n-1})\\
&\qquad=\sum_{\substack{x_1+\cdots+x_n=z_n\\x_t\geq0}}
\prod_{t=1}^{n}p_t^{*z_{n-t}}(x_t).
\end{aligned}
\end{equation}
Consequently,
\begin{equation}\label{eq:joint-pmf}
\begin{aligned}
&\Pr(Z_1=z_1,\ldots,Z_n=z_n)\\
&=\prod_{r=1}^{n}\left[
\sum_{\substack{x_1+\cdots+x_r=z_r\\x_t\geq0}}
\prod_{t=1}^{r}p_t^{*z_{r-t}}(x_t)
\right].
\end{aligned}
\end{equation}
\end{proposition}
\begin{proof}
Conditionally on the previous cohort sizes, the age-$t$ contribution to $Z_n$ is a sum of $z_{n-t}$ independent copies of $X(t)$, with distribution $p_t^{*z_{n-t}}$. These contributions are conditionally independent over $t$ by Assumption~\ref{ass:ind}. Convolving their distributions gives \eqref{eq:kernel}. The joint probability formula follows from the chain rule.
\end{proof}

\begin{remark}
The kernel in \eqref{eq:kernel} depends on the complete history through the vector $(z_0,\ldots,z_{n-1})$. Thus the scalar process $(Z_n)$ is generally not a time-homogeneous first-order Markov chain. A finite-dimensional Markov representation becomes available when the maximum reproductive age is bounded; see Section~\ref{sec:asymptotic}.
\end{remark}

\section{Martingale representation and second-order structure}\label{sec:second}
Assume in this section that $\E[X(t)^2]<\infty$ for each age under consideration. Define the innovation
\begin{equation}\label{eq:innovation}
\varepsilon_n=Z_n-\sum_{t=1}^{n}m_tZ_{n-t},\qquad n\geq1,
\end{equation}
and put $d_n=\E[\varepsilon_n^2]$.

\begin{lemma}[Innovation variance and orthogonality]\label{lem:innovation}
For every $n\geq1$,
\begin{equation}\label{eq:d-rec}
d_n=\sum_{t=1}^{n}\sigma_t^2a_{n-t}.
\end{equation}
Moreover, $\E[\varepsilon_n\varepsilon_j]=0$ whenever $j<n$.
\end{lemma}
\begin{proof}
Conditionally on $\cF_{n-1}$, the innovation is the sum of centered independent age-specific offspring contributions. Hence
\[
\E[\varepsilon_n\mid\cF_{n-1}]=0,
\qquad
\E[\varepsilon_n^2\mid\cF_{n-1}]
=\sum_{t=1}^{n}\sigma_t^2Z_{n-t}.
\]
Taking expectations gives \eqref{eq:d-rec}. If $j<n$, then $\varepsilon_j$ is $\cF_{n-1}$-measurable, and
\[
\E[\varepsilon_n\varepsilon_j]
=\E\left[\varepsilon_j\E(\varepsilon_n\mid\cF_{n-1})\right]=0.
\]
\end{proof}

\begin{theorem}[Linear innovation representation]\label{thm:linear-rep}
For every $n\geq1$,
\begin{equation}\label{eq:linear-rep}
Z_n-a_n=\sum_{j=1}^{n}a_{n-j}\varepsilon_j.
\end{equation}
\end{theorem}
\begin{proof}
The definition of $\varepsilon_n$ gives the random renewal equation
\[
Z_n=\sum_{t=1}^{n}m_tZ_{n-t}+\varepsilon_n.
\]
Subtract the deterministic renewal equation \eqref{eq:mean-rec}. The resulting triangular recursion for $Z_n-a_n$ has the same renewal coefficients $m_t$ and forcing term $\varepsilon_n$. Iterating it, the coefficient of $\varepsilon_j$ at time $n$ is the renewal sequence coefficient $a_{n-j}$, with $a_0=1$. This proves \eqref{eq:linear-rep}.
\end{proof}

\begin{theorem}[Covariance and variance formulae]\label{thm:covariance}
For $n,r\geq1$,
\begin{equation}\label{eq:covariance}
\Cov(Z_n,Z_r)=
\sum_{j=1}^{\min(n,r)}a_{n-j}a_{r-j}d_j,
\end{equation}
where $d_j$ is given by \eqref{eq:d-rec}. In particular,
\begin{equation}\label{eq:variance}
\Var(Z_n)=\sum_{j=1}^{n}a_{n-j}^{2}
\left(\sum_{t=1}^{j}\sigma_t^2a_{j-t}\right).
\end{equation}
All pairwise covariances between birth-cohort sizes are nonnegative.
\end{theorem}
\begin{proof}
Centering \eqref{eq:linear-rep} gives
\[
Z_n-a_n=\sum_{j=1}^{n}a_{n-j}\varepsilon_j.
\]
Multiply the analogous representations for $Z_n-a_n$ and $Z_r-a_r$, take expectations, and use the orthogonality in Lemma~\ref{lem:innovation}. Only terms with equal innovation indices remain, yielding \eqref{eq:covariance}. Set $r=n$ and substitute \eqref{eq:d-rec} to obtain \eqref{eq:variance}. Since $a_k\geq0$ and $d_j\geq0$, the covariance is nonnegative.
\end{proof}

\begin{corollary}[Correlation formula]\label{cor:correlation}
Whenever both variances are positive,
\[
\operatorname{Corr}(Z_n,Z_r)=
\frac{\displaystyle\sum_{j=1}^{\min(n,r)}a_{n-j}a_{r-j}d_j}
{\displaystyle\left(\sum_{j=1}^{n}a_{n-j}^2d_j\right)^{1/2}
\left(\sum_{j=1}^{r}a_{r-j}^2d_j\right)^{1/2}}.
\]
\end{corollary}

\begin{remark}[Cross-cohort terms]
A direct variance calculation from \eqref{eq:population} must account for dependence among the random cohort sizes $Z_{n-t}$. The variance is not, in general, the sum of the separate age-contribution variances. Formula \eqref{eq:variance} incorporates the cross-cohort dependence through the common innovations.
\end{remark}

\section{Arbitrary-order mixed moments and cumulants}\label{sec:allmom}

\subsection{Mixed factorial moments}
For integers $x,q\geq0$, let
\[
\fall{x}{q}=x(x-1)\cdots(x-q+1),\qquad \fall{x}{0}=1.
\]
For each age $t$, define the offspring factorial moments
\[
\mu_{t,q}=\E[\fall{X(t)}{q}],\qquad \mu_{t,0}=1.
\]
For $n\geq1$ and $\boldsymbol r=(r_1,\ldots,r_n)\in\Nzero^n$, write
\begin{equation}\label{eq:F-def}
F_n(\boldsymbol r)=\E\left[\prod_{i=1}^{n}\fall{Z_i}{r_i}\right],
\qquad F_0=1.
\end{equation}

\begin{definition}[Offspring-sum coefficients]\label{def:B}
For age $t\geq1$, integer $q\geq0$, and integer $x\geq0$, define
\begin{equation}\label{eq:B-def}
B_{t,q}(x)=q!\,[u^q]
\left(1+\sum_{r=1}^{q}\frac{\mu_{t,r}}{r!}u^r\right)^x,
\end{equation}
where $[u^q]$ denotes coefficient extraction. For $q=0$, set $B_{t,0}(x)=1$.
\end{definition}

\begin{lemma}[Factorial moments of an offspring sum]\label{lem:B}
If $X_1(t),\ldots,X_x(t)$ are independent copies of $X(t)$, then
\[
B_{t,q}(x)=\E\left[\fall{X_1(t)+\cdots+X_x(t)}{q}\right].
\]
Moreover,
\begin{equation}\label{eq:B-expand}
B_{t,q}(x)=\sum_{k=0}^{q}b_{t,q,k}\fall{x}{k},
\qquad
b_{t,q,k}=\frac{q!}{k!}[u^q]
\left(\sum_{r=1}^{q}\frac{\mu_{t,r}}{r!}u^r\right)^k.
\end{equation}
\end{lemma}
\begin{proof}
For a nonnegative integer-valued random variable $X(t)$, the binomial identity gives, at the level of coefficients through degree $q$,
\[
\E[(1+u)^{X(t)}]
=\sum_{r=0}^{q}\frac{\mu_{t,r}}{r!}u^r
\]
when terms of degree greater than $q$ are discarded. Independence raises this polynomial to the power $x$, and multiplication of the coefficient of $u^q$ by $q!$ gives the first identity. For the second, put
\[
A_t(u)=\sum_{r=1}^{q}\frac{\mu_{t,r}}{r!}u^r
\]
and use
\[
(1+A_t(u))^x=\sum_{k=0}^{x}\frac{\fall{x}{k}}{k!}A_t(u)^k.
\]
Extracting the coefficient of $u^q$ yields \eqref{eq:B-expand}. Since $A_t(u)^k$ has no terms of degree below $k$, coefficients with $k>q$ vanish.
\end{proof}

For illustration, the first three polynomials are
\[
\begin{aligned}
B_{t,1}(x)&=x\mu_{t,1},\\
B_{t,2}(x)&=x\mu_{t,2}+x(x-1)\mu_{t,1}^2,\\
B_{t,3}(x)&=x\mu_{t,3}+3x(x-1)\mu_{t,1}\mu_{t,2}
+x(x-1)(x-2)\mu_{t,1}^3.
\end{aligned}
\]

\begin{theorem}[Finite recursion for arbitrary mixed factorial moments]\label{thm:factorial-moments}
Fix $n\geq1$ and a multi-index $\boldsymbol r=(r_1,\ldots,r_n)\in\Nzero^n$, and put $p=|\boldsymbol r|=\sum_{i=1}^n r_i$. Suppose that $\mu_{t,q}<\infty$ for $1\leq t\leq n$ and $0\leq q\leq p$. For $n\geq2$, define, in each summand,
\[
\boldsymbol r'=(r_1+k_1-h_1,\ldots,r_{n-1}+k_{n-1}-h_{n-1}).
\]
Then
\begin{equation}\label{eq:factorial-recursion}
\begin{aligned}
F_n(\boldsymbol r)
={}&\sum_{\substack{q_1+\cdots+q_n=r_n\\q_t\geq0}}
\frac{r_n!}{q_1!\cdots q_n!}\,\mu_{n,q_n}\\
&\quad\times\sum_{k_1=0}^{q_{n-1}}\cdots\sum_{k_{n-1}=0}^{q_1}
\left(\prod_{i=1}^{n-1}b_{n-i,q_{n-i},k_i}\right)\\
&\quad\times\sum_{h_1=0}^{\min(r_1,k_1)}\cdots
\sum_{h_{n-1}=0}^{\min(r_{n-1},k_{n-1})}
\left[\prod_{i=1}^{n-1}
\binom{r_i}{h_i}\binom{k_i}{h_i}h_i!\right]
F_{n-1}(\boldsymbol r').
\end{aligned}
\end{equation}
For $n=1$, $F_1((r_1))=\mu_{1,r_1}$.
\end{theorem}
\begin{proof}
Write
\[
Y_{n,t}=\sum_{j=1}^{Z_{n-t}}X_{n-t,j}(t),
\qquad Z_n=\sum_{t=1}^{n}Y_{n,t}.
\]
Conditionally on $\cF_{n-1}$, the variables $Y_{n,1},\ldots,Y_{n,n}$ are independent. The factorial multinomial identity gives
\[
\fall{Z_n}{r_n}
=\sum_{q_1+\cdots+q_n=r_n}
\frac{r_n!}{q_1!\cdots q_n!}
\prod_{t=1}^{n}\fall{Y_{n,t}}{q_t}.
\]
Taking conditional expectations and applying Lemma~\ref{lem:B},
\[
\E\left[\fall{Z_n}{r_n}\mid\cF_{n-1}\right]
=\sum_{q_1+\cdots+q_n=r_n}
\frac{r_n!}{q_1!\cdots q_n!}
\prod_{t=1}^{n}B_{t,q_t}(Z_{n-t}).
\]
Multiply by $\prod_{i=1}^{n-1}\fall{Z_i}{r_i}$ and take expectations. Since $Z_0=1$, the factor corresponding to $t=n$ is
\[
B_{n,q_n}(1)=\mu_{n,q_n}.
\]
For a previous cohort $Z_i$, the relevant parental age at calendar time $n$ is $n-i$. Consequently, its offspring-sum polynomial is $B_{n-i,q_{n-i}}(Z_i)$, and the expansion in \eqref{eq:B-expand} reads
\[
B_{n-i,q_{n-i}}(Z_i)
=\sum_{k_i=0}^{q_{n-i}}
 b_{n-i,q_{n-i},k_i}\fall{Z_i}{k_i}.
\]
This explains the reversed age index $q_{n-i}$ in the inner sum: $q_t$ is indexed by parental age $t$, whereas $k_i$ is indexed by the earlier cohort $Z_i$. Finally, use
\begin{equation}\label{eq:falling-product}
\fall{x}{a}\fall{x}{b}
=\sum_{h=0}^{\min(a,b)}
\binom{a}{h}\binom{b}{h}h!\fall{x}{a+b-h},
\end{equation}
which follows by classifying pairs of ordered selections of sizes $a$ and $b$ according to the number $h$ of common selected elements. Applying \eqref{eq:falling-product} to each product $\fall{Z_i}{r_i}\fall{Z_i}{k_i}$ gives exactly \eqref{eq:factorial-recursion}.
\end{proof}

\begin{remark}[Finite-order closure]\label{rem:order-closure}
Every mixed factorial moment on the right-hand side of \eqref{eq:factorial-recursion} has total order at most $p$. Indeed,
\[
\begin{aligned}
\sum_{i=1}^{n-1}(r_i+k_i-h_i)
&\leq \sum_{i=1}^{n-1}r_i+\sum_{i=1}^{n-1}k_i\\
&\leq \sum_{i=1}^{n-1}r_i+\sum_{i=1}^{n-1}q_{n-i}
\leq p.
\end{aligned}
\]
Thus moments of total order at most $p$ at time $n$ can be computed from offspring factorial moments through order $p$ and earlier mixed moments of order at most $p$. The recursion is finite for every fixed $n$ and $p$.
\end{remark}

\begin{corollary}[Ordinary mixed moments]\label{cor:ordinary}
Let $S(a,r)$ denote a Stirling number of the second kind, defined by
\[
x^a=\sum_{r=0}^{a}S(a,r)\fall{x}{r}.
\]
Then for nonnegative integers $a_1,\ldots,a_n$,
\begin{equation}\label{eq:ordinary}
\E\left[\prod_{i=1}^{n}Z_i^{a_i}\right]
=\sum_{r_1=0}^{a_1}\cdots\sum_{r_n=0}^{a_n}
\left(\prod_{i=1}^{n}S(a_i,r_i)\right)
F_n(\boldsymbol r),
\end{equation}
where $\boldsymbol r=(r_1,\ldots,r_n)$ in each summand.
\end{corollary}
\begin{proof}
Apply the Stirling expansion separately to each factor $Z_i^{a_i}$, multiply the finite sums, and take expectations.
\end{proof}

\subsection{A recursion for mixed cumulants}
We next derive a finite recursion for ordinary mixed cumulants. The argument is formulated with truncated formal power series. Thus no exponential-moment assumption is imposed: finite ordinary moments of the required order are sufficient.

Fix an integer $R\geq1$. For a multi-index $\boldsymbol r=(r_1,\ldots,r_n)$, write $|\boldsymbol r|=\sum_i r_i$, $\boldsymbol r!=\prod_i r_i!$, and $\boldsymbol\theta^{\boldsymbol r}=\prod_i\theta_i^{r_i}$. Let $\mathbb{R}[[\theta_1,\ldots,\theta_n]]_{\leq R}$ denote the ring of formal power series in which all monomials of total degree greater than $R$ are discarded. Define
\[
\mathfrak M_n^{[R]}(\boldsymbol\theta)
=\sum_{|\boldsymbol r|\leq R}
\E\left[\prod_{i=1}^{n}Z_i^{r_i}\right]
\frac{\boldsymbol\theta^{\boldsymbol r}}{\boldsymbol r!},
\qquad
\mathfrak K_n^{[R]}=\log\mathfrak M_n^{[R]},
\]
where $\mathfrak M_n^{[R]}$ and its logarithm are both interpreted in this truncated ring. For $|\boldsymbol r|\leq R$, set
\begin{equation}\label{eq:C-def}
C_n(\boldsymbol r)
=\boldsymbol r!\,[\boldsymbol\theta^{\boldsymbol r}]
\mathfrak K_n^{[R]}(\boldsymbol\theta).
\end{equation}
These coefficients are the ordinary joint cumulants of the cohort variables, with repeated arguments allowed. For example,
\[
C_2((1,1))=\Cov(Z_1,Z_2),\qquad
C_2((2,1))=\cum(Z_1,Z_1,Z_2).
\]
For each age $t$, let $c_{t,q}$ denote the $q$th ordinary cumulant of $X(t)$, and define the truncated offspring cumulant series
\[
\psi_t^{[R]}(z)=\sum_{q=1}^{R}c_{t,q}\frac{z^q}{q!}.
\]
All substitutions, products, and coefficient extractions below are made in the appropriate truncated formal power-series ring. In particular, $\psi_t^{[R]}(0)=0$.

\begin{theorem}[Mixed-cumulant recursion]\label{thm:cumulants}
Let $n\geq2$, $q\geq1$, and $\boldsymbol\alpha=(\alpha_1,\ldots,\alpha_{n-1})\in\Nzero^{n-1}$. Put $R=|\boldsymbol\alpha|+q$. Suppose that
\[
\E[X(t)^R]<\infty,\qquad 1\leq t\leq n.
\]
Then
\begin{equation}\label{eq:cumulant-recursion}
\begin{aligned}
C_n((\boldsymbol\alpha,q))
={}&\1_{\{\boldsymbol\alpha=\boldsymbol0\}}c_{n,q}\\
&+q!
\sum_{\substack{\boldsymbol\beta\in\Nzero^{n-1}\\
1\leq|\boldsymbol\beta|\leq q}}
\frac{C_{n-1}(\boldsymbol\alpha+\boldsymbol\beta)}
{\boldsymbol\beta!}
[z^q]\prod_{i=1}^{n-1}
\bigl(\psi_{n-i}^{[R]}(z)\bigr)^{\beta_i}.
\end{aligned}
\end{equation}
For $n=1$, $C_1((q))=c_{1,q}$. The recursion is finite for every fixed $q$ and involves only cumulants of total order at most $R$.
\end{theorem}
\begin{proof}
The conditional pgf identity and the independence assumption imply the corresponding formal moment-generating identity through total degree $R$:
\[
\mathfrak M_n^{[R]}(\boldsymbol\theta,z)
=\exp\bigl(\psi_n^{[R]}(z)\bigr)
\mathfrak M_{n-1}^{[R]}\bigl(\theta_1+\psi_{n-1}^{[R]}(z),\ldots,
\theta_{n-1}+\psi_1^{[R]}(z)\bigr),
\]
where the identity is read modulo all monomials of total degree greater than $R$. Taking the formal logarithm gives
\begin{equation}\label{eq:K-recursion}
\begin{aligned}
\mathfrak K_n^{[R]}(\theta_1,\ldots,\theta_{n-1},z)
={}&\psi_n^{[R]}(z)\\
&+\mathfrak K_{n-1}^{[R]}\bigl(\theta_1+\psi_{n-1}^{[R]}(z),\ldots,
\theta_{n-1}+\psi_1^{[R]}(z)\bigr).
\end{aligned}
\end{equation}

By the definition of the formal cumulants,
\[
\mathfrak K_{n-1}^{[R]}(\boldsymbol u)
=\sum_{|\boldsymbol k|\leq R}
C_{n-1}(\boldsymbol k)\frac{\boldsymbol u^{\boldsymbol k}}
{\boldsymbol k!}.
\]
Set $u_i=\theta_i+\psi_{n-i}^{[R]}(z)$. Differentiating formally with respect to $\theta_i$ to order $\alpha_i$ for each $i=1,\ldots,n-1$, and then setting $\boldsymbol\theta=\boldsymbol0$, yields
\[
\begin{aligned}
&\left.\partial_{\theta_1}^{\alpha_1}\cdots
\partial_{\theta_{n-1}}^{\alpha_{n-1}}
\mathfrak K_n^{[R]}(\boldsymbol\theta,z)\right|_{\boldsymbol\theta=\boldsymbol0}\\
&=\1_{\{\boldsymbol\alpha=\boldsymbol0\}}\psi_n^{[R]}(z)
+\sum_{\substack{\boldsymbol\beta\in\Nzero^{n-1}\\
|\boldsymbol\alpha+\boldsymbol\beta|\leq R}}
\frac{C_{n-1}(\boldsymbol\alpha+\boldsymbol\beta)}
{\boldsymbol\beta!}
\prod_{i=1}^{n-1}
\bigl(\psi_{n-i}^{[R]}(z)\bigr)^{\beta_i}.
\end{aligned}
\]
Each $\psi_{n-i}^{[R]}(z)$ has zero constant term. Hence the product in the summand has order at least $|\boldsymbol\beta|$ in $z$. Its coefficient of $z^q$ is therefore zero if $|\boldsymbol\beta|>q$. The term $\boldsymbol\beta=\boldsymbol0$ is constant in $z$ and contributes nothing because $q\geq1$. Extracting the coefficient of $z^q$ and multiplying by $q!$ proves \eqref{eq:cumulant-recursion}, since
\[
q![z^q]\psi_n^{[R]}(z)=c_{n,q}.
\]
All coefficients used in this calculation have total degree at most $R$; thus moments of the offspring variables through order $R$ are sufficient. This completes the proof.
\end{proof}

For clarity, the first two cases of the recursion are recorded explicitly. Let $e_i$ denote the $i$th coordinate vector in $\mathbb R^{n-1}$. For $q=1$,
\[
C_n((\boldsymbol\alpha,1))
=\1_{\{\boldsymbol\alpha=\boldsymbol0\}}m_n
+\sum_{i=1}^{n-1}m_{n-i}
C_{n-1}(\boldsymbol\alpha+e_i).
\]
For $q=2$,
\[
\begin{aligned}
C_n((\boldsymbol\alpha,2))
={}&\1_{\{\boldsymbol\alpha=\boldsymbol0\}}\sigma_n^2
+\sum_{i=1}^{n-1}\sigma_{n-i}^2
C_{n-1}(\boldsymbol\alpha+e_i)\\
&+\sum_{i=1}^{n-1}m_{n-i}^2
C_{n-1}(\boldsymbol\alpha+2e_i)\\
&+2\sum_{1\leq i<j\leq n-1}m_{n-i}m_{n-j}
C_{n-1}(\boldsymbol\alpha+e_i+e_j).
\end{aligned}
\]
The $q=2$ identity is consistent with the second-order decomposition, including all cross-covariance contributions.

\section{Finite memory and supercritical limits}\label{sec:asymptotic}

\subsection{Finite-memory state representation}
If there is an integer $L\geq1$ such that $X(t)=0$ almost surely for every $t>L$, then the vector
\[
W_n=(Z_n,Z_{n-1},\ldots,Z_{n-L+1})^\top
\]
is a time-homogeneous first-order Markov chain. Indeed, only the most recent $L$ birth cohorts can contribute to the next cohort, and the remaining coordinates are shifted deterministically. The scalar process $(Z_n)$ is generally not first-order Markov.

\subsection{A finite-age Perron limit theorem}
Throughout this subsection, assume that $X(t)=0$ almost surely for $t>L$. Write
\[
m_t=\E[X(t)],\qquad \sigma_t^2=\Var(X(t)),\qquad 1\leq t\leq L,
\]
and set $Z_k=0$ for $k<0$. Define the finite-age mean matrix
\begin{equation}\label{eq:finiteL-A}
A=\begin{pmatrix}
m_1&m_2&\cdots&m_{L-1}&m_L\\
1&0&\cdots&0&0\\
0&1&\cdots&0&0\\
\vdots&\ddots&\ddots&\vdots&\vdots\\
0&\cdots&0&1&0
\end{pmatrix},
\qquad e_1=(1,0,\ldots,0)^\top.
\end{equation}
The innovation is now
\begin{equation}\label{eq:finiteL-innovation}
\varepsilon_n=Z_n-\sum_{t=1}^{L}m_tZ_{n-t},
\end{equation}
and the state recursion becomes
\begin{equation}\label{eq:finiteL-state-rec}
W_n=AW_{n-1}+e_1\varepsilon_n,
\qquad W_0=e_1.
\end{equation}
Conditionally on $\mathcal F_{n-1}$, the innovation is a sum of centered independent offspring variables. Consequently,
\begin{equation}\label{eq:finiteL-eps-cond}
\E[\varepsilon_n\mid\mathcal F_{n-1}]=0,
\qquad
\E[\varepsilon_n^2\mid\mathcal F_{n-1}]
=\sum_{t=1}^{L}\sigma_t^2Z_{n-t}.
\end{equation}

\begin{assumption}[Primitive supercritical finite-age case]\label{ass:finiteL}
Assume $m_L>0$ and
\[
\gcd\{t\in\{1,\ldots,L\}:m_t>0\}=1.
\]
Let $\lambda=\rho(A)$ be the Perron--Frobenius eigenvalue and assume $\lambda>1$.
\end{assumption}

\begin{lemma}[Primitivity of the companion mean matrix]\label{lem:finiteL-primitive}
The matrix $A$ in \eqref{eq:finiteL-A} is primitive if and only if $m_L>0$ and the greatest common divisor of the set $\{t:m_t>0\}$ is one.
\end{lemma}
\begin{proof}
Consider the directed graph associated with $A$, with an edge from $i$ to $j$ whenever $A_{ij}>0$. The subdiagonal entries give the edges $i\to i-1$ for $i=2,\ldots,L$. If $m_L>0$, the edge $1\to L$ is also present. These edges make the graph strongly connected: every vertex reaches vertex $1$ by successive subdiagonal edges, and vertex $1$ reaches every vertex by first moving to $L$ and then following the subdiagonal edges. Thus $A$ is irreducible. Every directed cycle can be decomposed into cycles that start at vertex $1$, use an edge $1\to t$ corresponding to $m_t>0$, and return to vertex $1$ along the subdiagonal chain. Such a basic cycle has length $t$. Hence the period of the irreducible graph is exactly $\gcd\{t:m_t>0\}$. An irreducible nonnegative matrix is primitive precisely when its period is one. Conversely, if $m_L=0$, the graph is not strongly connected on all $L$ coordinates, and if the displayed gcd exceeds one, the graph is periodic. This proves the equivalence.
\end{proof}

Under Assumption~\ref{ass:finiteL}, the Perron--Frobenius theorem gives a simple eigenvalue $\lambda>1$ whose modulus strictly exceeds that of every other eigenvalue. Let $r$ and $\ell$ be the corresponding positive right and left eigenvectors, normalized by
\begin{equation}\label{eq:finiteL-eigen-norm}
r_1=1,\qquad \ell^\top r=1.
\end{equation}
The companion structure gives the explicit right eigenvector
\begin{equation}\label{eq:finiteL-r}
r=(1,\lambda^{-1},\ldots,\lambda^{-(L-1)})^\top,
\end{equation}
and the characteristic equation is
\begin{equation}\label{eq:finiteL-characteristic}
\lambda^L=m_1\lambda^{L-1}+m_2\lambda^{L-2}+\cdots+m_L,
\qquad
\sum_{t=1}^{L}m_t\lambda^{-t}=1.
\end{equation}
In particular, $\ell_1>0$. Define
\begin{equation}\label{eq:finiteL-perron-mart}
M_n=\lambda^{-n}\ell^\top W_n.
\end{equation}

\begin{lemma}[Finite-age Perron martingale]\label{lem:finiteL-mart}
Under Assumption~\ref{ass:finiteL} and finite first offspring moments, $(M_n,\mathcal F_n)$ is a nonnegative martingale and
\begin{equation}\label{eq:finiteL-mart-increment}
M_n-M_{n-1}=\ell_1\lambda^{-n}\varepsilon_n.
\end{equation}
Moreover, $\E[M_n]=\ell_1$ for all $n\geq0$.
\end{lemma}
\begin{proof}
Taking conditional expectations in \eqref{eq:finiteL-state-rec} and using \eqref{eq:finiteL-eps-cond},
\[
\E[W_n\mid\mathcal F_{n-1}]=AW_{n-1}.
\]
It follows that
\[
\E[M_n\mid\mathcal F_{n-1}]
=\lambda^{-n}\ell^\top AW_{n-1}
=\lambda^{-(n-1)}\ell^\top W_{n-1}=M_{n-1}.
\]
Also, by \eqref{eq:finiteL-state-rec},
\[
M_n=\lambda^{-n}\ell^\top AW_{n-1}
+\lambda^{-n}\ell^\top e_1\varepsilon_n
=M_{n-1}+\ell_1\lambda^{-n}\varepsilon_n.
\]
The entries of $\ell$ and $W_n$ are nonnegative, so $M_n\geq0$. Finally, $W_0=e_1$ and hence $M_0=\ell_1$; the martingale property gives the expectation identity.
\end{proof}

\begin{lemma}[Mean growth and innovation variance]\label{lem:finiteL-mean-growth}
Under Assumption~\ref{ass:finiteL},
\begin{equation}\label{eq:finiteL-mean-vector}
\lambda^{-n}\E[W_n]=\lambda^{-n}A^ne_1\longrightarrow \ell_1r.
\end{equation}
In particular,
\begin{equation}\label{eq:finiteL-mean-asymptotic}
a_n:=\E[Z_n]\sim\ell_1\lambda^n.
\end{equation}
If the offspring variables have finite second moments, then
\begin{equation}\label{eq:finiteL-dn}
d_n:=\E[\varepsilon_n^2]
=\sum_{t=1}^{\min(n,L)}\sigma_t^2a_{n-t}
=O(\lambda^n).
\end{equation}
\end{lemma}
\begin{proof}
The mean recursion gives $\E[W_n]=A^ne_1$. Since $A$ is primitive, its Perron eigenvalue is simple and strictly dominant, and the Perron projection is $P=r\ell^\top$. Thus $\lambda^{-n}A^n\to P$, which proves \eqref{eq:finiteL-mean-vector} and \eqref{eq:finiteL-mean-asymptotic}. Taking expectations in the conditional variance identity \eqref{eq:finiteL-eps-cond} gives the expression for $d_n$. Since $L$ is finite and $a_k=O(\lambda^k)$, the bound $d_n=O(\lambda^n)$ follows.
\end{proof}

\begin{lemma}[Moment growth for the finite-age process]\label{lem:finiteL-moment-growth}
Assume Assumption~\ref{ass:finiteL}. Fix an integer $p\geq2$ and suppose
\[
\E[X(t)^p]<\infty,\qquad t=1,\ldots,L.
\]
Then
\begin{equation}\label{eq:finiteL-moment-growth}
\E[Z_n^p]=O(\lambda^{np}),\qquad
\|\varepsilon_n\|_p=O(\lambda^{n/2}).
\end{equation}
\end{lemma}
\begin{proof}
We use induction on the integer moment order. The first moment bound follows from Lemma~\ref{lem:finiteL-mean-growth}. Fix $p\geq2$ and assume that the population moment bound has been established for every positive integer strictly smaller than $p$. Put $q=\lceil p/2\rceil$ and
\[
S_n=\sum_{t=1}^{L}Z_{n-t}.
\]
For $p=2$, $q=1$ and the required bound for $S_n^q$ follows from the mean estimate. For $p>2$, one has $q<p$, so the induction hypothesis and the finite value of $L$ give
\[
\E[S_n^q]=O(\lambda^{nq}).
\]
By monotonicity of $L^r$ norms,
\[
\E[S_n^{p/2}]
\leq\bigl(\E[S_n^q]\bigr)^{p/(2q)}
=O(\lambda^{np/2}),
\qquad
\E[S_n]=O(\lambda^n).
\]
Conditionally on $\mathcal F_{n-1}$, $\varepsilon_n$ is a sum of $S_n$ independent centered random variables, each with one of the finitely many age-specific centered offspring laws. Rosenthal's inequality~\cite{Rosenthal1970} (and the conditional variance identity when $p=2$) therefore yields a constant $C_p$ such that
\[
\begin{aligned}
\E[|\varepsilon_n|^p\mid\mathcal F_{n-1}]
&\leq C_p\left[
\left(\sum_{t=1}^{L}\sigma_t^2Z_{n-t}\right)^{p/2}
+\sum_{t=1}^{L}\E|X(t)-m_t|^p Z_{n-t}
\right]\\
&\leq C_p'(S_n^{p/2}+S_n).
\end{aligned}
\]
Taking expectations gives $\E|\varepsilon_n|^p=O(\lambda^{np/2})$, or equivalently $\|\varepsilon_n\|_p=O(\lambda^{n/2})$.

It remains to close the induction for $Z_n$. Iterating \eqref{eq:finiteL-state-rec},
\[
W_n=A^nW_0+\sum_{j=1}^{n}A^{n-j}e_1\varepsilon_j.
\]
The Perron--Frobenius spectral gap implies that $\|A^k\|\leq C\lambda^k$ for all $k\geq0$. Minkowski's inequality gives
\[
\begin{aligned}
\|W_n\|_p
&\leq C\lambda^n\|W_0\|_p
+C\sum_{j=1}^{n}\lambda^{n-j}\|\varepsilon_j\|_p\\
&\leq C\lambda^n
+C\sum_{j=1}^{n}\lambda^{n-j}\lambda^{j/2}
=O(\lambda^n),
\end{aligned}
\]
because $\lambda>1$. Thus $\E[Z_n^p]=O(\lambda^{np})$, completing the induction.
\end{proof}

\begin{theorem}[Exact $L^p$ criterion and finite-age state limit]\label{thm:finiteL-Lp}
Assume Assumption~\ref{ass:ind} and Assumption~\ref{ass:finiteL}, and fix an integer $p\geq2$. Then the following are equivalent:
\begin{enumerate}[label=(\roman*)]
\item $\E[X(t)^p]<\infty$ for every $t=1,\ldots,L$;
\item the Perron martingale $(M_n)$ converges in $L^p$;
\item the normalized state vector $\lambda^{-n}W_n$ converges in $L^p$.
\end{enumerate}
When these conditions hold, there exists $M_\infty\in L^p$ such that
\begin{equation}\label{eq:finiteL-vector-Lp}
\lambda^{-n}W_n\xrightarrow{L^p}rM_\infty,
\qquad
M_n\xrightarrow{L^p}M_\infty,
\qquad \E[M_\infty]=\ell_1.
\end{equation}
Consequently, with $U=M_\infty/\ell_1$,
\begin{equation}\label{eq:finiteL-normalized-Lp}
\frac{Z_n}{a_n}\xrightarrow{L^p}U,
\qquad \E[U]=1,
\end{equation}
and, for every fixed $h\geq0$,
\begin{equation}\label{eq:finiteL-fixed-lag-Lp}
\frac{Z_{n-h}}{a_{n-h}}\xrightarrow{L^p}U.
\end{equation}
\end{theorem}
\begin{proof}
We first prove that (i) implies (ii). By Lemma~\ref{lem:finiteL-moment-growth},
\[
\|\varepsilon_n\|_p=O(\lambda^{n/2}).
\]
The increment formula in Lemma~\ref{lem:finiteL-mart} gives
\[
\|M_n-M_{n-1}\|_p
=\ell_1\lambda^{-n}\|\varepsilon_n\|_p
=O(\lambda^{-n/2}).
\]
The series of these norms is summable. Hence $(M_n)$ is Cauchy in $L^p$ and converges in $L^p$ to some $M_\infty$. Since convergence in $L^p$ implies convergence in $L^1$, the expectations converge and $\E[M_\infty]=\ell_1$.

We next prove that (ii) implies (i). Suppose that $M_n$ converges in $L^p$. In particular, it converges in $L^1$, and the martingale convergence identity gives
\[
\E[M_\infty\mid\mathcal F_t]=M_t,
\qquad t\geq0.
\]
Let
\[
X_t=X_{0,1}(t),\qquad t=1,\ldots,L,
\]
be the age-specific offspring counts of the initial ancestor. For each fixed $t$, $X_t$ is $\mathcal F_t$-measurable and $Z_t\geq X_t$. All entries of $\ell$ and all coordinates of $W_t$ are nonnegative. Therefore,
\[
M_t=\lambda^{-t}\ell^\top W_t
\geq\lambda^{-t}\ell_1Z_t
\geq\lambda^{-t}\ell_1X_t.
\]
Conditional Jensen's inequality now yields
\[
\E[M_\infty^p]
\geq\E\left[\bigl(\E[M_\infty\mid\mathcal F_t]\bigr)^p\right]
=\E[M_t^p]
\geq(\lambda^{-t}\ell_1)^p\E[X_t^p].
\]
The left-hand side is finite, so $\E[X_t^p]<\infty$. Since $X_t$ has the distribution of $X(t)$, this holds for every $t=1,\ldots,L$. Thus (ii) implies (i).

We now prove that (i) implies (iii). Let
\[
P=r\ell^\top,\qquad Q=I-P,\qquad B=A-\lambda P.
\]
The normalization $\ell^\top r=1$ implies $P^2=P$, and the eigenvector identities imply $AP=PA=\lambda P$ and $PB=BP=0$. Thus
\[
A^n=\lambda^nP+B^n.
\]
All eigenvalues of $B$ have modulus strictly less than $\lambda$. Choose $\theta<\lambda$ larger than their moduli. By finite-dimensional spectral theory, there is a constant $C$ such that
\begin{equation}\label{eq:finiteL-B-bound}
\|B^n\|\leq C\theta^n,\qquad n\geq0.
\end{equation}
Applying $Q$ to the iterated state recursion gives
\[
QW_n=B^nQW_0+\sum_{j=1}^{n}B^{n-j}Qe_1\varepsilon_j.
\]
Consequently, by Minkowski's inequality and Lemma~\ref{lem:finiteL-moment-growth},
\[
\begin{aligned}
\|\lambda^{-n}QW_n\|_p
&\leq C(\theta/\lambda)^n
+C\lambda^{-n}\sum_{j=1}^{n}\theta^{n-j}\lambda^{j/2}\\
&=C(\theta/\lambda)^n
+C\lambda^{-n/2}\sum_{k=0}^{n-1}(\theta/\sqrt\lambda)^k
\longrightarrow0.
\end{aligned}
\]
The last convergence holds whether $\theta/\sqrt\lambda$ is smaller than, equal to, or larger than one, since $\theta<\lambda$. On the Perron eigenspace,
\[
\lambda^{-n}PW_n=r\lambda^{-n}\ell^\top W_n=rM_n
\xrightarrow{L^p}rM_\infty.
\]
Combining the two components proves
\[
\lambda^{-n}W_n\xrightarrow{L^p}rM_\infty,
\]
so (i) implies (iii).

Finally, (iii) implies (ii) because
\[
M_n=\ell^\top(\lambda^{-n}W_n)
\]
and multiplication by the fixed row vector $\ell^\top$ is continuous in $L^p$. We have proved all three implications and hence the equivalence.

The first coordinate of the state-vector convergence and Lemma~\ref{lem:finiteL-mean-growth} give
\[
\frac{Z_n}{a_n}
=\frac{\lambda^{-n}Z_n}{\lambda^{-n}a_n}
\xrightarrow{L^p}\frac{M_\infty}{\ell_1}=U.
\]
Replacing $n$ by $n-h$ proves the fixed-lag assertion. The mean of $U$ equals one because $\E[M_\infty]=\ell_1$. The proof is complete.
\end{proof}

\begin{remark}[The finite-type branching representation]\label{rem:finiteL-multitype}
The finite-memory state recursion is also a finite-type branching representation. A type-$i$ individual, $i=1,\ldots,L$, produces $X(i)$ type-1 individuals and, if $i<L$, one type-$(i+1)$ individual representing its continuation to the next age. With row vectors of type counts, the mean matrix is $A^\top$. Thus the Perron martingale mechanism is consistent with the classical theory of indecomposable multitype Galton--Watson processes; see Kesten and Stigum~\cite{KestenStigum1966}. The specific statements below exploit the companion-matrix structure to express the limits directly in terms of birth cohorts and age-specific offspring moments.
\end{remark}

\subsection{\texorpdfstring{Mixed cumulants of fixed-lag cohorts for arbitrary finite $L$}{Mixed cumulants of fixed-lag cohorts for arbitrary finite L}}
We now use the vector limit to obtain explicit asymptotic statements for the dependence between several birth cohorts. For a random variable $Y$ with a finite $r$th moment, write $\kappa_r(Y)=\cum_r(Y)$.

\begin{theorem}[Fixed-lag mixed-cumulant limit for finite $L$]\label{thm:finiteL-cumulant-limit}
Assume the hypotheses of Theorem~\ref{thm:finiteL-Lp} for an integer $r\geq2$. Let $h_1,\ldots,h_r$ be fixed nonnegative integers. Then
\begin{equation}\label{eq:finiteL-cumulant-limit}
\lim_{n\to\infty}
\frac{\cum(Z_{n-h_1},\ldots,Z_{n-h_r})}
{\prod_{j=1}^{r}a_{n-h_j}}
=\kappa_r(U),
\end{equation}
where $U=M_\infty/\ell_1$ is the mean-one limit in \eqref{eq:finiteL-normalized-Lp}. The limit is independent of the particular fixed lag vector $(h_1,\ldots,h_r)$.
\end{theorem}
\begin{proof}
For each fixed $h\geq0$, Theorem~\ref{thm:finiteL-Lp} gives
\[
Y_{n,h}:=Z_{n-h}/a_{n-h}\xrightarrow{L^r}U.
\]
Let $B$ be a nonempty subset of $\{1,\ldots,r\}$. We claim that
\[
\E\left[\prod_{j\in B}Y_{n,h_j}\right]
\longrightarrow \E[U^{|B|}].
\]
To verify this, write the difference of the products as a telescoping sum, replacing one factor at a time by $U$. Each summand contains one difference converging to zero in $L^r$ and at most $r-1$ other factors whose $L^r$ norms are uniformly bounded. H\"older's inequality, with exponents chosen so that the product is integrable, gives convergence of each summand in $L^1$. (If $|B|<r$, the same conclusion follows by the continuous embedding $L^r\subseteq L^{|B|}$ on a probability space.)

For any $r$ random variables with finite moments through order $r$, their joint cumulant is the finite partition sum
\begin{equation}\label{eq:finiteL-cumulant-partition}
\cum(Y_1,\ldots,Y_r)
=\sum_{\pi\in\mathcal P_r}
(-1)^{|\pi|-1}(|\pi|-1)!
\prod_{B\in\pi}\E\left[\prod_{j\in B}Y_j\right],
\end{equation}
where $\mathcal P_r$ denotes the set of partitions of $\{1,\ldots,r\}$. The preceding product-moment convergence allows passage to the limit in every term of this finite sum. Hence
\[
\cum(Y_{n,h_1},\ldots,Y_{n,h_r})
\longrightarrow\cum(\underbrace{U,\ldots,U}_{r\text{ arguments}})
=\kappa_r(U).
\]
By multilinearity of cumulants, the left-hand side is exactly the normalized expression in \eqref{eq:finiteL-cumulant-limit}. This proves the theorem.
\end{proof}

\subsection{A smoothing transform and explicit cumulants of the limit}
The fixed-lag theorem identifies the limit of mixed cohort cumulants in terms of the cumulants of $U$. The next result characterizes those cumulants directly from the age-specific offspring laws.

\begin{theorem}[Smoothing-transform equation for the Perron limit]\label{thm:finiteL-smoothing}
Under the assumptions of Theorem~\ref{thm:finiteL-Lp} with $p=2$, let $X_t=X_{0,1}(t)$ be the offspring count of the initial ancestor at age $t$. There exist mutually independent copies $U_{t,j}$ of $U$, independent of $(X_1,\ldots,X_L)$, such that
\begin{equation}\label{eq:finiteL-smoothing}
U\ \overset{d}=\sum_{t=1}^{L}\lambda^{-t}
\sum_{j=1}^{X_t}U_{t,j}.
\end{equation}
The sums are finite almost surely, and $\E[U]=1$.
\end{theorem}
\begin{proof}
Construct the branching tree by first exposing the age-specific offspring counts $X_t$ of the initial ancestor. Each child born at time $t$ initiates an independent copy of the original process. Denote the cohort process in the subtree rooted at child $j$ born at time $t$ by $Z_k^{(t,j)}$, with $Z_0^{(t,j)}=1$. The branching property and the disjointness of descendant subtrees give, for every $n\geq L$,
\begin{equation}\label{eq:finiteL-root-decomp}
Z_n=\sum_{t=1}^{L}\sum_{j=1}^{X_t}Z_{n-t}^{(t,j)}.
\end{equation}
For each subtree, use the same Perron projection $P=r\ell^\top$ and $Q=I-P$ as in the proof of Theorem~\ref{thm:finiteL-Lp}. Define its first-coordinate spectral remainder by
\[
R_k^{(t,j)}=e_1^\top\lambda^{-k}QW_k^{(t,j)}.
\]
Since $r_1=1$ and $M_k^{(t,j)}=\lambda^{-k}\ell^\top W_k^{(t,j)}$, the exact decomposition is
\[
\lambda^{-k}Z_k^{(t,j)}=M_k^{(t,j)}+R_k^{(t,j)}.
\]
The estimate in the proof of Theorem~\ref{thm:finiteL-Lp} gives
\[
\|R_k^{(t,j)}\|_2
\leq C(\theta/\lambda)^k
+C\lambda^{-k/2}\sum_{u=0}^{k-1}(\theta/\sqrt{\lambda})^u.
\]
Because $\theta<\lambda$, there are constants $C'<\infty$ and $0<\rho<1$ such that the right-hand side is at most $C'\rho^k$ (choose $\rho$ strictly between $\max\{\theta/\lambda,\lambda^{-1/2}\}$ and $1$). Consequently,
\[
\sum_{k\geq0}\mathbb E|R_k^{(t,j)}|
\leq\sum_{k\geq0}\|R_k^{(t,j)}\|_2
\leq C'\sum_{k\geq0}\rho^k<\infty.
\]
Tonelli's theorem implies $\sum_k|R_k^{(t,j)}|<\infty$ almost surely, and hence $R_k^{(t,j)}\to0$ almost surely. The martingale $M_k^{(t,j)}$ converges almost surely as well, so the first-coordinate limit holds almost surely for every fixed subtree. The limits are independent copies of $M_\infty$ and are independent of the root offspring vector. Since each $X_t$ is finite almost surely and there are only finitely many ages, we may pass to the limit in the finite sum in \eqref{eq:finiteL-root-decomp} to obtain
\[
M_\infty\overset{d}=\sum_{t=1}^{L}\lambda^{-t}
\sum_{j=1}^{X_t}M_\infty^{(t,j)}.
\]
Dividing by $\ell_1$ and setting $U_{t,j}=M_\infty^{(t,j)}/\ell_1$ proves \eqref{eq:finiteL-smoothing}. The mean-one property follows from Theorem~\ref{thm:finiteL-Lp}.
\end{proof}

For each $t$, let $\kappa_{t,k}=\cum_k(X(t))$ denote the $k$th ordinary cumulant of the offspring count at age $t$, and put
\[
\gamma_r=\cum_r(U),\qquad \gamma_1=1.
\]
The smoothing transform yields a finite triangular recursion for $(\gamma_r)_{r\geq2}$.

\begin{theorem}[Recursive formula for the cumulants of the limit]\label{thm:finiteL-limit-cumulants}
Fix an integer $r\geq2$ and assume $\E[X(t)^r]<\infty$ for every $t=1,\ldots,L$. Define
\[
D_r=\sum_{t=1}^{L}m_t\lambda^{-tr}.
\]
Then $D_r<1$ and
\begin{equation}\label{eq:finiteL-limit-cumulant-rec}
\gamma_r=
\frac{r!}{1-D_r}
\sum_{t=1}^{L}\sum_{k=2}^{r}
\frac{\kappa_{t,k}}{k!}
[z^r]
\left(\sum_{j=1}^{r-1}
\frac{\gamma_j\lambda^{-tj}}{j!}z^j\right)^k.
\end{equation}
In particular,
\begin{equation}\label{eq:finiteL-gamma2}
\gamma_2=
\frac{\displaystyle\sum_{t=1}^{L}\sigma_t^2\lambda^{-2t}}
{\displaystyle1-\sum_{t=1}^{L}m_t\lambda^{-2t}}.
\end{equation}
For the third cumulant,
\begin{equation}\label{eq:finiteL-gamma3}
\gamma_3=
\frac{\displaystyle\sum_{t=1}^{L}\lambda^{-3t}
\left(3\sigma_t^2\gamma_2+\kappa_{t,3}\right)}
{\displaystyle1-\sum_{t=1}^{L}m_t\lambda^{-3t}}.
\end{equation}
For completeness, the fourth cumulant is
\begin{equation}\label{eq:finiteL-gamma4}
\gamma_4=
\frac{\displaystyle\sum_{t=1}^{L}\lambda^{-4t}
\left(\kappa_{t,4}+6\kappa_{t,3}\gamma_2
+4\sigma_t^2\gamma_3+3\sigma_t^2\gamma_2^2\right)}
{\displaystyle1-\sum_{t=1}^{L}m_t\lambda^{-4t}}.
\end{equation}
\end{theorem}
\begin{proof}
By \eqref{eq:finiteL-characteristic},
\[
D_1=\sum_{t=1}^{L}m_t\lambda^{-t}=1.
\]
Since $\lambda>1$, $m_t\geq0$, and at least one $m_t$ is positive, for every $r\geq2$,
\[
0\leq D_r=\sum_{t=1}^{L}m_t\lambda^{-tr}
<\sum_{t=1}^{L}m_t\lambda^{-t}=1.
\]
Thus the denominator in \eqref{eq:finiteL-limit-cumulant-rec} is strictly positive.

We derive the recursion using formal cumulant series, so no moment-generating function is assumed to exist in a neighborhood of zero. Work in the ring of formal power series truncated after degree $r$. Let
\[
K_U(z)=\sum_{j=1}^{r}\gamma_j\frac{z^j}{j!},
\qquad
K_t(z)=\sum_{k=1}^{r}\kappa_{t,k}\frac{z^k}{k!}.
\]
The independence in the smoothing equation implies, at the level of formal cumulant series through degree $r$,
\begin{equation}\label{eq:finiteL-formal-smoothing}
K_U(z)=\sum_{t=1}^{L}K_t\bigl(K_U(\lambda^{-t}z)\bigr).
\end{equation}
Indeed, conditionally on $X_t$, the moment series of the age-$t$ sum is the $X_t$th power of the moment series of $\lambda^{-t}U$; taking the formal logarithm first over the independent copies and then over $X_t$ gives the composition $K_t(K_U(\lambda^{-t}z))$. The age-specific sums are independent across $t$, so their cumulant series add.

Expanding the right-hand side of \eqref{eq:finiteL-formal-smoothing},
\[
K_t(K_U(\lambda^{-t}z))
=\sum_{k=1}^{r}\frac{\kappa_{t,k}}{k!}
\left(\sum_{j=1}^{r}
\frac{\gamma_j\lambda^{-tj}}{j!}z^j\right)^k
\quad\pmod{z^{r+1}}.
\]
The coefficient of $z^r$ in the term $k=1$ is
\[
\frac{m_t\gamma_r\lambda^{-tr}}{r!}.
\]
For $k\geq2$, the coefficient of $z^r$ cannot involve $\gamma_r$, because every factor in the $k$th power has degree at least one, and using the degree-$r$ term from one factor would force the total degree to exceed $r$. Thus, for $k\geq2$, $\gamma_r$ may be omitted from the inner sum and the coefficient is unchanged if the sum is truncated at $j=r-1$. Equating the coefficients of $z^r$ on both sides of \eqref{eq:finiteL-formal-smoothing}, multiplying by $r!$, and moving the $k=1$ terms to the left gives
\[
(1-D_r)\gamma_r
=r!\sum_{t=1}^{L}\sum_{k=2}^{r}
\frac{\kappa_{t,k}}{k!}
[z^r]
\left(\sum_{j=1}^{r-1}
\frac{\gamma_j\lambda^{-tj}}{j!}z^j\right)^k,
\]
which is \eqref{eq:finiteL-limit-cumulant-rec}.

For $r=2$, only $k=2$ contributes. Since $\gamma_1=1$ and $\kappa_{t,2}=\sigma_t^2$,
\[
2!\frac{\kappa_{t,2}}{2!}
[z^2](\lambda^{-t}z)^2
=\sigma_t^2\lambda^{-2t},
\]
which proves \eqref{eq:finiteL-gamma2}. For $r=3$, the $k=2$ term contributes $3\sigma_t^2\gamma_2\lambda^{-3t}$ and the $k=3$ term contributes $\kappa_{t,3}\lambda^{-3t}$. This gives \eqref{eq:finiteL-gamma3}. For $r=4$, the coefficients of $z^4$ in the terms $k=2,3,4$ give, respectively, the contributions
\[
\lambda^{-4t}\left(4\sigma_t^2\gamma_3+3\sigma_t^2\gamma_2^2\right),
\qquad
6\lambda^{-4t}\kappa_{t,3}\gamma_2,
\qquad
\lambda^{-4t}\kappa_{t,4},
\]
which proves \eqref{eq:finiteL-gamma4}.
\end{proof}

\begin{corollary}[Explicit asymptotic covariance and fixed-lag correlation]\label{cor:finiteL-covariance}
Assume the hypotheses of Theorem~\ref{thm:finiteL-Lp} with $p=2$. For any fixed nonnegative integers $h,k$,
\begin{equation}\label{eq:finiteL-cov-limit}
\frac{\Cov(Z_{n-h},Z_{n-k})}
{a_{n-h}a_{n-k}}
\longrightarrow
\frac{\displaystyle\sum_{t=1}^{L}\sigma_t^2\lambda^{-2t}}
{\displaystyle1-\sum_{t=1}^{L}m_t\lambda^{-2t}}
=\gamma_2.
\end{equation}
Equivalently, using $a_n\sim\ell_1\lambda^n$,
\begin{equation}\label{eq:finiteL-raw-cov-limit}
\lambda^{-2n+h+k}\Cov(Z_{n-h},Z_{n-k})
\longrightarrow\ell_1^2\gamma_2.
\end{equation}
If at least one age-specific offspring distribution is non-degenerate, then $\gamma_2>0$ and
\begin{equation}\label{eq:finiteL-corr-limit}
\operatorname{Corr}(Z_{n-h},Z_{n-k})\longrightarrow1.
\end{equation}
If all age-specific offspring variables are deterministic, then $U=1$ almost surely and all limiting cumulants of order at least two vanish; in that degenerate case the correlations need not be defined.
\end{corollary}
\begin{proof}
Theorem~\ref{thm:finiteL-cumulant-limit} with $r=2$ gives the normalized covariance limit, and Theorem~\ref{thm:finiteL-limit-cumulants} gives the explicit value $\gamma_2$. For completeness, the same value follows directly from the martingale increments. Since they are orthogonal,
\[
\Var(M_\infty)
=\ell_1^2\sum_{n=1}^{\infty}\lambda^{-2n}d_n
=\ell_1^2\sum_{t=1}^{L}\sigma_t^2\lambda^{-2t}
\sum_{j=0}^{\infty}a_j\lambda^{-2j}.
\]
The mean renewal equation has generating function
\[
\sum_{j=0}^{\infty}a_jz^j
=\frac{1}{1-\sum_{t=1}^{L}m_tz^t}
\]
for $z$ in a neighborhood of zero and, by monotone convergence, at $z=\lambda^{-2}$, since $\sum_t m_t\lambda^{-2t}<1$. Dividing by $\ell_1^2$ yields the displayed expression for $\Var(U)=\gamma_2$.

If at least one $\sigma_t^2$ is positive, the numerator is positive and the denominator is positive, so $\gamma_2>0$. Applying \eqref{eq:finiteL-cov-limit} to the covariance and to each of the two variances gives the correlation limit. If all offspring variables are deterministic, the entire process and its normalized limit are deterministic, so $U=1$ and its higher cumulants vanish.
\end{proof}

\begin{remark}[Scope of fixed-lag limits]\label{rem:finiteL-fixed-lags}
The limits in Theorem~\ref{thm:finiteL-cumulant-limit} and Corollary~\ref{cor:finiteL-covariance} hold for fixed lags. Theorem~\ref{thm:growing-gap-mse} below treats a different regime in which both the gap and the earlier observation time diverge, and gives a quantitative mean-square asymptotic under a simple-root spectral expansion.
\end{remark}

\subsection{Mean-square asymptotics for growing cohort gaps}
We now quantify the distance between two mean-normalized cohorts when both their observation times and their separation increase. Write
\[
D_{n,h}=\frac{Z_n}{a_n}-\frac{Z_{n-h}}{a_{n-h}},\qquad m=n-h.
\]
The exact innovation representation gives a useful decomposition into innovations shared by the two cohorts and innovations appearing only in the later cohort.

For this subsection, put
\[
d_j=\mathbb E[\varepsilon_j^2],\qquad
S=\sum_{t=1}^{L}\sigma_t^2\lambda^{-t}.
\]
Under the primitive supercritical assumptions, $a_j\sim c\lambda^j$ with $c=\ell_1$, and hence
\begin{equation}\label{eq:dn-asymptotic-growing}
\frac{d_j}{\lambda^j}\longrightarrow cS.
\end{equation}
We impose an additional spectral assumption on the scalar renewal sequence. It is satisfied, for example, when the dominant non-Perron eigenvalue of the companion matrix is real, simple, and has a nonzero coefficient in $e_1^\top A^ne_1$.

\begin{assumption}[Simple real subdominant root]\label{ass:simple-subroot}
There exist $b\ne0$, a real number $\beta$ with $0<|\beta|<\lambda$, and a number $\eta$ with $0<\eta<|\beta|$, such that
\begin{equation}\label{eq:mean-two-root-expansion}
a_k=c\lambda^k+b\beta^k+O(\eta^k),\qquad k\to\infty,
\end{equation}
where $c=\ell_1>0$. The remainder bound is understood with a constant independent of $k$.
\end{assumption}

Set $q=\beta/\lambda$ and $f_k=a_k/(c\lambda^k)$. The following uniform estimate is the technical point needed to pass to limits in the sums below.

\begin{lemma}[Uniform ratio estimates for the renewal sequence]\label{lem:uniform-renewal-ratios}
Under Assumption~\ref{ass:simple-subroot}, there are constants $C<\infty$ and $m_0$ such that, for all $m\geq m_0$, all $h\geq1$, and all $1\leq j\leq m$,
\begin{equation}\label{eq:uniform-ratio-bound}
\left|\frac{f_{m+h-j}}{f_{m+h}}-\frac{f_{m-j}}{f_m}\right|
\leq C|q|^{m-j}.
\end{equation}
Moreover, for each fixed $j\geq1$, along any sequence $m\to\infty$, $h\to\infty$,
\begin{equation}\label{eq:uniform-ratio-pointwise}
q^{-m}\left(\frac{f_{m+h-j}}{f_{m+h}}-\frac{f_{m-j}}{f_m}\right)
\longrightarrow -\frac{b}{c}q^{-j}(1-q^j).
\end{equation}
\end{lemma}
\begin{proof}
Let $\delta=b/c$. Dividing \eqref{eq:mean-two-root-expansion} by $c\lambda^k$ gives
\[
f_k=1+\delta q^k+R_k,\qquad |R_k|\leq C_0\xi^k,
\]
where $\xi=\eta/\lambda<|q|$. Enlarging constants to cover the finitely many initial indices, we also have
\begin{equation}\label{eq:f-uniform-decay}
|f_k-1|\leq C_1|q|^k,\qquad k\geq0.
\end{equation}
Since $f_k\to1$, there are $m_0$ and $c_0>0$ such that $f_k\geq c_0$ for all $k\geq m_0$. In particular, the denominators $f_m$ and $f_{m+h}$ are bounded away from zero for $m\geq m_0$.

For $k\geq m_0$ and $1\leq j\leq k$, it follows from \eqref{eq:f-uniform-decay} that
\[
\begin{aligned}
\left|\frac{f_{k-j}}{f_k}-1\right|
&=\frac{|f_{k-j}-f_k|}{f_k}\\
&\leq c_0^{-1}\bigl(|f_{k-j}-1|+|f_k-1|\bigr)\\
&\leq C_2\bigl(|q|^{k-j}+|q|^k\bigr)
\leq C_3|q|^{k-j}.
\end{aligned}
\]
Apply this inequality first with $k=m+h$ and then with $k=m$. The triangle inequality gives
\[
\left|\frac{f_{m+h-j}}{f_{m+h}}-\frac{f_{m-j}}{f_m}\right|
\leq C_3\bigl(|q|^{m+h-j}+|q|^{m-j}\bigr)
\leq 2C_3|q|^{m-j},
\]
which proves \eqref{eq:uniform-ratio-bound}.

For the pointwise limit, fix $j$. From the expansion of $f_k$ and the identity
\[
\frac{1+x}{1+y}=1+(x-y)+O(x^2+y^2+|xy|)
\]
for $x,y\to0$, we obtain
\begin{equation}\label{eq:fixed-j-ratio-expansion}
\frac{f_{k-j}}{f_k}
=1+\delta q^{k-j}(1-q^j)+O(\zeta^{k-j}),
\end{equation}
where one may choose $\zeta<|q|$, for example any number strictly larger than $\max\{\xi,|q|^2\}$. Subtract \eqref{eq:fixed-j-ratio-expansion} at $k=m$ from the same expansion at $k=m+h$ and divide by $q^m$. The leading term becomes
\[
\delta q^{-j}(1-q^j)(q^h-1),
\]
which converges to $-\delta q^{-j}(1-q^j)$ because $h\to\infty$. The two remainder terms, after division by $|q|^m$, tend to zero: for the term with $k=m$, this follows from $(\zeta/|q|)^{m-j}\to0$; for $k=m+h$, it follows from
\[
|q|^{-m}\zeta^{m+h-j}
=|q|^{h-j}\left(\frac{\zeta}{|q|}\right)^{m+h-j}\longrightarrow0.
\]
This proves \eqref{eq:uniform-ratio-pointwise}.
\end{proof}

\begin{theorem}[Three-regime spectral asymptotics for growing gaps]\label{thm:growing-gap-mse}
Assume the primitive supercritical finite-age conditions, finite second offspring moments, $S>0$, and Assumption~\ref{ass:simple-subroot}. Let $m=n-h$ and suppose
\[
m\longrightarrow\infty,\qquad h\longrightarrow\infty.
\]
Then
\begin{equation}\label{eq:growing-gap-three-regime}
\mathbb E[D_{n,h}^2]\sim
\begin{cases}
C_<\lambda^{-m},&\beta^2<\lambda,\\[1mm]
C_=\,m\lambda^{-m},&\beta^2=\lambda,\\[1mm]
C_>\left(\dfrac{|\beta|}{\lambda}\right)^{2m},&\beta^2>\lambda,
\end{cases}
\end{equation}
where
\begin{align}
C_<&=cS\left\{\frac{1}{\lambda-1}
+\sum_{r=0}^{\infty}\lambda^r(1-f_r)^2\right\},\label{eq:C-under}\\
C_=&=\frac{Sb^2}{c},\label{eq:C-critical}\\
C_>&=\frac{b^2}{c^2}\sum_{j=1}^{\infty}
\frac{d_j}{|\beta|^{2j}}(1-q^j)^2.
\label{eq:C-over}
\end{align}
All three constants are finite and strictly positive. The series in $C_<$ converges when $\beta^2<\lambda$, and the series in $C_>$ converges when $\beta^2>\lambda$.
\end{theorem}
\begin{proof}
The linear innovation representation and orthogonality of the innovations give the exact identity
\begin{align}
\mathbb E[D_{n,h}^2]
={}&\sum_{j=1}^{m}d_j
\left(\frac{a_{n-j}}{a_n}-\frac{a_{m-j}}{a_m}\right)^2
+\sum_{j=m+1}^{n}d_j\left(\frac{a_{n-j}}{a_n}\right)^2\\
=:&\ T_{1,m,h}+T_{2,m,h}.
\label{eq:growing-gap-exact}
\end{align}
The second term is always of order at most $\lambda^{-m}$. Indeed, for all $k\geq0$ there are constants $C_1,C_2$ such that $a_k\leq C_1\lambda^k$, while $a_{m+h}\geq C_2\lambda^{m+h}$ for all sufficiently large $m+h$, and $d_{m+s}\leq C_1\lambda^{m+s}$. Writing $s=j-m$ gives
\[
T_{2,m,h}
=\sum_{s=1}^{h}d_{m+s}\left(\frac{a_{h-s}}{a_{m+h}}\right)^2
\leq C\sum_{s=1}^{h}\lambda^{-m-s}
=O(\lambda^{-m}).
\label{eq:T2-bound}
\]

For the first term, use $r=m-j$ and the definition of $f_k$. Then
\begin{equation}\label{eq:T1-r}
T_{1,m,h}
=\lambda^{-2m}\sum_{r=0}^{m-1}d_{m-r}\lambda^{2r}
\left(\frac{f_{h+r}}{f_{m+h}}-\frac{f_r}{f_m}\right)^2.
\end{equation}
Put
\[
G_{m,h,r}=\frac{f_{h+r}}{f_{m+h}}-\frac{f_r}{f_m}.
\]
The same estimate as in the first part of Lemma~\ref{lem:uniform-renewal-ratios}, now with the indices $h+r$ and $m+h$, together with \eqref{eq:f-uniform-decay}, yields
\begin{equation}\label{eq:G-r-bound}
|G_{m,h,r}|\leq C|q|^r,
\qquad 0\leq r\leq m-1,
\end{equation}
with $C$ independent of $m,h,r$ for all sufficiently large $m$. For each fixed $r$,
\[
G_{m,h,r}\longrightarrow 1-f_r,
\qquad
\frac{d_{m-r}}{\lambda^{m-r}}\longrightarrow cS.
\]

Suppose first that $\beta^2<\lambda$. Multiplying \eqref{eq:T1-r} by $\lambda^m$ gives
\[
\lambda^mT_{1,m,h}
=\sum_{r=0}^{m-1}
\frac{d_{m-r}}{\lambda^{m-r}}\lambda^rG_{m,h,r}^2.
\]
By \eqref{eq:G-r-bound}, each summand is bounded by
\[
C\lambda^r|q|^{2r}
=C\left(\frac{\beta^2}{\lambda}\right)^r,
\]
a summable sequence. Dominated convergence for series therefore gives
\[
\lambda^mT_{1,m,h}
\longrightarrow cS\sum_{r=0}^{\infty}\lambda^r(1-f_r)^2.
\]
Combining this with \eqref{eq:T2-bound} and
$\lambda^mT_{2,m,h}\to cS/(\lambda-1)$ requires a sharper limit for $T_2$, which follows by writing
\[
\lambda^mT_{2,m,h}
=\sum_{s=1}^{h}
\frac{d_{m+s}}{\lambda^{m+s}}
\lambda^{-s}
\left(\frac{f_{h-s}}{f_{m+h}}\right)^2.
\]
For each fixed $s$, the summand tends to $cS\lambda^{-s}$. It is bounded by $C\lambda^{-s}$ uniformly in $m,h,s$ because the $f_k$ are bounded and $f_{m+h}$ is bounded away from zero. Dominated convergence gives
\[
\lambda^mT_{2,m,h}\longrightarrow cS\sum_{s=1}^{\infty}\lambda^{-s}
=\frac{cS}{\lambda-1}.
\]
This proves the first regime.

Suppose next that $\beta^2=\lambda$, so $|q|=\lambda^{-1/2}$. The summands in the expression for $\lambda^mT_{1,m,h}$ are uniformly bounded by a constant, by \eqref{eq:G-r-bound}. We show that they converge uniformly to $Sb^2/c$ away from the two ends of the summation interval. Specifically, for every $\epsilon>0$ there is $R$ such that, whenever $r\geq R$, $m-r\geq R$, and $h\geq R$,
\[
\left|\frac{d_{m-r}}{\lambda^{m-r}}-cS\right|<\epsilon,
\qquad
\left|\frac{G_{m,h,r}}{q^r}+\frac bc\right|<\epsilon.
\]
The first inequality follows from \eqref{eq:dn-asymptotic-growing}. For the second, write $\delta=b/c$ and use
$f_k=1+\delta q^k+R_k$, where $|R_k|\leq C\xi^k$ for some $\xi<|q|$. The exact identities
\[
\begin{aligned}
\frac{f_{h+r}}{f_{m+h}}-1
&=\frac{\delta(q^{h+r}-q^{m+h})+R_{h+r}-R_{m+h}}{f_{m+h}},\\
\frac{f_r}{f_m}-1
&=\frac{\delta(q^r-q^m)+R_r-R_m}{f_m}.
\end{aligned}
\]
show, after subtraction and division by $q^r$, that
\[
\left|\frac{G_{m,h,r}}{q^r}+\delta\right|
\leq C\left(|q|^h+|q|^{m-r}+|q|^m
 +(\xi/|q|)^r+\xi^h+\xi^{m-r}\right)
\]
whenever $r\geq R$ and $m-r,h\geq R$. The constant is independent of $m,h,r$: this follows because $f_m$ and $f_{m+h}$ are bounded away from zero, while $|q|<1$ and $\xi/|q|<1$. The right-hand side tends to zero uniformly on the indicated middle range as $R\to\infty$. Hence $G_{m,h,r}/q^r$ converges uniformly there to $-\delta=-b/c$. Consequently,
\[
\frac{d_{m-r}}{\lambda^{m-r}}\lambda^rG_{m,h,r}^2
\longrightarrow cS\frac{b^2}{c^2}=\frac{Sb^2}{c}
\]
uniformly on the same middle range.

More explicitly, in the critical case $|q|^2=\lambda^{-1}$, and \eqref{eq:G-r-bound} gives the uniform estimate
\[
0\leq
\frac{d_{m-r}}{\lambda^{m-r}}\lambda^rG_{m,h,r}^2
\leq C\lambda^r|q|^{2r}=C,
\qquad 0\leq r\leq m-1,
\]
where $C$ is independent of $m,h,r$. Thus the first and last $R$ terms contribute at most $2CR=O(R)$, uniformly in $m,h$. After division by $m$, these endpoint contributions vanish for each fixed $R$. Letting first $m\to\infty$ and then $R\to\infty$ yields
\[
\frac{\lambda^mT_{1,m,h}}{m}\longrightarrow\frac{Sb^2}{c}.
\]
Since $T_{2,m,h}=O(\lambda^{-m})$, the second regime follows.

Finally, suppose that $\beta^2>\lambda$. Re-index $T_{1,m,h}$ by $j=m-r$:
\[
T_{1,m,h}
=\sum_{j=1}^{m}d_j\lambda^{-2j}
\left(\frac{f_{m+h-j}}{f_{m+h}}-\frac{f_{m-j}}{f_m}\right)^2.
\]
For each fixed $j$, Lemma~\ref{lem:uniform-renewal-ratios} gives
\[
\left(\frac{\lambda}{|\beta|}\right)^{2m}
\lambda^{-2j}
\left(\frac{f_{m+h-j}}{f_{m+h}}-\frac{f_{m-j}}{f_m}\right)^2
\longrightarrow
\frac{b^2}{c^2}\frac{(1-q^j)^2}{|\beta|^{2j}}.
\]
Moreover, the uniform bound \eqref{eq:uniform-ratio-bound} implies that the left-hand side is at most
\[
C\lambda^{-2j}|q|^{-2j}=C|\beta|^{-2j}.
\]
Thus the summand after multiplication by $d_j$ is dominated by
$C d_j|\beta|^{-2j}$. This is summable because $d_j=O(\lambda^j)$ and $\lambda/\beta^2<1$. Extending the summands by zero for $j>m$, dominated convergence for series yields
\[
\left(\frac{\lambda}{|\beta|}\right)^{2m}T_{1,m,h}
\longrightarrow
\frac{b^2}{c^2}\sum_{j=1}^{\infty}
\frac{d_j}{|\beta|^{2j}}(1-q^j)^2=C_>.
\]
Since $\lambda^{-m}=o((|\beta|/\lambda)^{2m})$ in this regime, the bound on $T_2$ shows that it is negligible. To verify strict positivity of the constant, note that $S>0$ implies that there is an age $t_0\in\{1,\ldots,L\}$ with $\sigma_{t_0}^2>0$. Since $a_0=1$,
\[
d_{t_0}=\sum_{t=1}^{t_0}\sigma_t^2a_{t_0-t}
\geq\sigma_{t_0}^2a_0>0.
\]
Moreover, $0<|q|<1$ implies $q^{t_0}\ne1$. Hence the $j=t_0$ summand in \eqref{eq:C-over} is strictly positive, and therefore $C_> >0$. This proves the third regime and completes the proof.
\end{proof}

\begin{corollary}[Correlation for arbitrary diverging observation times]\label{cor:growing-gap-correlation}
Under the assumptions of Theorem~\ref{thm:growing-gap-mse},
\[
\operatorname{Corr}(Z_n,Z_{n-h})\longrightarrow1.
\]
In particular, this holds for $h\to\infty$ with $h=o(n)$, but the condition $h=o(n)$ is not required: it is enough that both $h\to\infty$ and $n-h\to\infty$.
\end{corollary}
\begin{proof}
Theorem~\ref{thm:finiteL-Lp} with $p=2$ implies that both $Z_n/a_n$ and $Z_{n-h}/a_{n-h}$ converge in $L^2$ to the same non-degenerate random variable $U$, since both indices tend to infinity. Hence their variances converge to $\operatorname{Var}(U)>0$ and their covariance converges to the same value. The correlation is unchanged by positive deterministic normalization, so the result follows. The three-regime theorem gives the sharper rate for the $L^2$ distance.
\end{proof}

The two-age case provides a transparent illustration of the spectral condition used above.

\subsection{A two-age illustration}
When $L=2$, the mean matrix is
\[
A=\begin{pmatrix}m_1&m_2\\1&0\end{pmatrix},
\]
and its eigenvalues are
\[
\lambda=\frac{m_1+\sqrt{m_1^2+4m_2}}2,
\qquad
\lambda_-=-\frac{m_2}{\lambda}.
\]
If $m_1,m_2>0$ and $m_1+m_2>1$, the matrix is primitive and supercritical. Its two eigenvalues are distinct, with $\lambda>0$ and $\beta=\lambda_-=-m_2/\lambda<0$. The scalar mean sequence has the exact representation
\[
 a_k=c\lambda^k+b\beta^k,\qquad
 c=\frac{m_1-\beta}{\lambda-\beta},\quad
 b=\frac{\lambda-m_1}{\lambda-\beta}.
\]
Both coefficients are nonzero when $m_2>0$, so Assumption~\ref{ass:simple-subroot} holds (the remainder is identically zero). This gives a directly checkable sufficient condition for the spectral assumption in dimension two.

The three regimes are all attainable within this two-age family. For instance, take $m_1=1$ and prescribe the Perron root $\lambda>1$ by setting $m_2=\lambda(\lambda-1)$. Then $\beta=1-\lambda$. The choices
\[
\lambda=2,\quad
\lambda=\frac{3+\sqrt5}{2},\quad
\lambda=3
\]
respectively give $\beta^2<\lambda$, $\beta^2=\lambda$, and $\beta^2>\lambda$. The corresponding values of $m_2$ are $2$, $2+\sqrt5$ (equivalently $\lambda(\lambda-1)$), and $6$. To realize each choice as a complete stochastic example satisfying the moment and non-degeneracy assumptions of Theorem~\ref{thm:growing-gap-mse}, take all age-specific offspring variables to be independent and let $X(1)\sim\operatorname{Poisson}(m_1)$ and $X(2)\sim\operatorname{Poisson}(m_2)$. Then all offspring moments are finite and
\[
S=\sigma_1^2\lambda^{-1}+\sigma_2^2\lambda^{-2}
=m_1\lambda^{-1}+m_2\lambda^{-2}=1,
\]
by the characteristic equation. Thus $S>0$ in all three cases. The threshold is not a vacuous possibility imposed by an abstract spectral assumption; it occurs for explicit positive age-specific means in the smallest nontrivial dimension. For $L\geq3$, non-Perron roots can be complex or have equal moduli, which is why the theorem states a simple-real-root hypothesis rather than silently treating all companion matrices as if they had a real subdominant root.

\begin{remark}[What changes for a complex dominant non-Perron root]\label{rem:complex-subroot}
The real-root assumption is substantive. Suppose instead that the dominant non-Perron roots of the real renewal sequence are a simple conjugate pair $\zeta,\overline\zeta$, with $\zeta=\rho e^{i\vartheta}$, $0<\vartheta<\pi$, and that
\[
 a_k=c\lambda^k+b\zeta^k+\overline b\,\overline\zeta^{\,k}+O(\eta^k),
 \qquad 0<\eta<\rho<\lambda.
\]
Then the normalized remainder is oscillatory rather than a single real geometric term. In the subcritical case $\rho^2<\lambda$, the series defining the analogue of $C_< $ still converges, because its summands are bounded by a constant times $(\rho^2/\lambda)^k$; the same dominated-convergence argument gives a positive constant after retaining the actual real sequence $f_k=a_k/(c\lambda^k)$. At the critical boundary $\rho^2=\lambda$, the leading summands contain a non-oscillatory component and oscillatory components. For a genuinely non-real pair, Ces\`aro averaging removes the latter, but the resulting constant depends on the conjugate-pair amplitude and must be derived with the corresponding cross terms. In the supercritical case $\rho^2>\lambda$, the normalized leading term generally retains a phase depending on $m\vartheta$; consequently a single constant asymptotic equivalent need not exist along the full sequence. One should instead expect an oscillatory leading factor, or subsequential equivalents when the phases converge. These cases are not covered by Theorem~\ref{thm:growing-gap-mse}; proving a unified complex-root theorem requires a separate formulation that accounts for phase interactions and, for tied roots or Jordan blocks, polynomial factors. We therefore state them as extensions rather than claiming that the real-root theorem already covers them.
\end{remark}

\section{Discussion and comparison with related limit theorems}\label{sec:discussion}
The results should be interpreted within the assumptions of the independent age-specific reproduction array. The conditional pgf factorization is a consequence of independence not only between distinct parents, but also between the different age-specific offspring counts of the same parent. A general Crump--Mode--Jagers reproduction point process allows dependence between births at different parental ages; the present formulae do not cover that more general dependence structure without replacing the product of age-wise pgfs by an appropriate joint reproduction transform.

The arbitrary-order recursions in Section~\ref{sec:allmom} are finite coefficient identities for the specified cohort model. Their role is to provide explicit recursive access to mixed factorial moments and mixed cumulants while retaining the dependence among cohorts. The general problem of deriving moment equations for age-dependent branching processes is much older and is treated, in broader settings, by renewal and point-process methods; see Mode~\cite{Mode1974}, Sevast'yanov and Chistyakov~\cite{SevastyanovChistyakov1971}, and Kimmel~\cite{Kimmel1983}. 
Fathi Vajargah and Moradi~\cite{FathiMoradi2012} explicitly studied a discrete-time period-dependent branching process with overlapping generations in which an individual may reproduce at several ages. Their earlier epidemic formulation~\cite{FathiMoradi2011} also made the age/period dependence of reproduction central to the model. The current manuscript is a substantial distributional and asymptotic development of that line of work. Fathi Vajargah and Moradi~\cite{FathiMoradi2012} in their paper focused on the mean, variance, marginal pgf calculations, and simulation, whereas the present paper develops the joint cohort law, the exact history-dependent transition kernel, arbitrary-order mixed factorial-moment and cumulant recursions, and explicit finite-age limit/cumulant formulae.

For orientation, the relationship with the most relevant strands of prior work is summarized in Table~\ref{tab:positioning}. The table is deliberately conservative: it identifies the principal emphasis of each source rather than asserting that a listed topic is absent from every part of the cited work.

\begin{table}[htbp]
\centering
\scriptsize
\setlength{\tabcolsep}{3pt}
\caption{Positioning of the present results relative to selected branching-process literature.}
\label{tab:positioning}
\begin{tabular}{>{\raggedright\arraybackslash}p{0.20\textwidth}>{\raggedright\arraybackslash}p{0.30\textwidth}>{\raggedright\arraybackslash}p{0.35\textwidth}}
\toprule
Work & Main setting or emphasis & Relation to the present manuscript\\
\midrule
Fathi Vajargah and Moradi (2012)~\cite{FathiMoradi2012} & Discrete-time period-dependent branching with overlapping generations; mean, variance, pgf calculations, and simulation & Direct model predecessor. The present paper develops joint cohort laws, exact transition probabilities, mixed moments/cumulants, and finite-age limits.\\
Mode (1974)~\cite{Mode1974} & Discrete-time generalized age-dependent branching; renewal equations for means and covariances & Establishes the broader age-dependent moment-equation context; the present formulas specialize to an independent-age repeated-reproduction array.\\
Kimmel (1983)~\cite{Kimmel1983} & Point-process methods for age- and time-dependent branching and finite-dimensional distributions & Provides a broad distributional framework; the present paper gives explicit coefficient and cohort recursions under stronger independence assumptions.\\
Kesten and Stigum (1966)~\cite{KestenStigum1966} & Perron limits for indecomposable multitype Galton--Watson processes & Supplies the classical finite-type Perron background; the present work translates the state limit into age-indexed cohort and offspring-moment formulas.\\
Kolesko and Sava-Huss (2023)~\cite{KoleskoSavaHuss2023} & Limit theorems for discrete multitype branching processes counted with characteristics & A recent general CMJ-related limit-theorem framework; the present results focus on fixed-lag cohort cumulants and the finite-age smoothing transform.\\
Present work & Independent-age repeated reproduction with finite maximum age $L$ & Explicit joint cohort pgf and kernel; finite coefficient recursions; exact $L^p$ criterion; common-limit cumulants; three-regime mean-square asymptotics for growing cohort gaps under a simple-root condition.\\
\bottomrule
\end{tabular}
\end{table}

For finite $L$, the state vector is a finite-type branching process with an age-progression mechanism. Thus the existence of a Perron martingale and its first-order direction are structurally related to the classical multitype Galton--Watson theory of Kesten and Stigum~\cite{KestenStigum1966}. Recent work on discrete multitype branching processes counted with characteristics also establishes law-of-large-numbers and central-limit results in a broader CMJ-related setting; see Kolesko and Sava-Huss~\cite{KoleskoSavaHuss2023}. Our finite-age limit theorem is not a replacement for these general results. Its purpose is to exploit the companion-matrix structure to obtain an exact age-by-age $L^p$ moment criterion, direct limits for birth cohorts, and a closed recursive expression for the cumulants of their common normalized limit. The cumulant recursion follows from the model's independent-age smoothing transform and is stated only under the finite-moment assumptions required at each order.

It is useful to compare the asymptotic results with two further developments.

First, Janson's paper on lattice Crump--Mode--Jagers processes was published online in 2019 and appears in volume 50 (2018) of \emph{Advances in Applied Probability}~\cite{Janson2019}. It studies second-order fluctuations of age distributions in a general supercritical lattice CMJ setting. In particular, the asymptotic behaviour depends on the location of the next characteristic root relative to the square-root scale of the Malthusian root, producing Gaussian or oscillatory fluctuation regimes under suitable conditions. Our result addresses a different object and a different scale: for a fixed finite maximum reproductive age, it gives an $L^p$ limit for the finite-dimensional birth-cohort state and a mean-square asymptotic for the difference between two mean-normalized cohorts whose observation times and separation both diverge. The three regimes in Theorem~\ref{thm:growing-gap-mse} are determined by the dominant non-Perron root of the scalar renewal sequence. We do not derive a central limit theorem for the population state, and our result should not be read as a replacement for Janson's fluctuation theory.

Second, Cheng's 2024 preprint on supercritical remaining-lifetime age-structured branching processes~\cite{Cheng2024}, subsequently published in 2026, studies measure-valued population states in which an individual's reproduction law depends on its remaining lifetime. It establishes Malthusian limit results for normalized random measures, gives an $L\log L$-type condition for almost-sure strengthening, and develops a central limit theorem with random normalization. The state variable, reproduction mechanism, and convergence questions differ from those here. In particular, the necessary-and-sufficient $L^p$ criterion proved in Theorem~\ref{thm:finiteL-Lp} concerns the Perron martingale generated by the finite-dimensional birth-cohort recursion; it is not a replacement for Cheng's condition for convergence of normalized random measures. The two results are complementary rather than directly comparable in strength.

The main limitations are also clear. Dependence between the age-specific offspring counts of a single parent is excluded. The finite-age Perron theorem requires a primitive mean matrix; periodic cases require a separate analysis along residue classes. The three-regime growing-gap theorem assumes that the scalar mean sequence has a simple real dominant non-Perron term with nonzero coefficient. It therefore does not cover complex conjugate roots, ties in the non-Perron spectral radius, or nontrivial Jordan blocks. Those cases can produce oscillatory factors or polynomial corrections and require separate statements. The finite-type representation also makes clear that the Perron direction itself is classical. The main additional asymptotic information here is the explicit mean-square distance between two increasingly separated normalized cohorts, including the transition at the square-root spectral threshold.

\section{Conclusion}\label{sec:conclusion}
We have studied the joint law and temporal dependence of birth cohorts in a discrete-time age-dependent branching model with repeated reproduction and independent age-specific offspring counts. The finite-time analysis gives a recursive joint pgf, an exact history-dependent transition kernel, and an orthogonal innovation representation. These yield covariance identities and finite coefficient recursions for arbitrary mixed factorial moments and cumulants.

For finite maximum reproductive age, the cohort state has a companion-matrix representation. In the primitive supercritical case, the Perron martingale converges in $L^p$ exactly when all age-specific offspring distributions have finite $p$th moments, for each integer $p\geq2$. The normalized cohorts share a common limit; its smoothing-transform equation gives a triangular recursion for all cumulants under the corresponding finite-moment assumptions.

The principal quantitative result concerns growing temporal gaps. Under a simple real subdominant-root expansion of the mean renewal sequence, we prove a three-regime equivalent for the mean-square distance between the normalized sizes of cohorts $n$ and $n-h$, when both $h$ and $n-h$ tend to infinity. The rate changes at the threshold $\beta^2=\lambda$: it is of order $\lambda^{-(n-h)}$ below the threshold, $(n-h)\lambda^{-(n-h)}$ at the threshold, and $(|\beta|/\lambda)^{2(n-h)}$ above it. The proof uses uniform ratio bounds and a justified passage of limits through the relevant series, including in the supercritical spectral regime.

The results rely on independence across parental ages and on a primitive finite-age mean matrix. The growing-gap theorem further assumes a simple real dominant non-Perron root in the scalar mean sequence. Extending the spectral analysis to complex roots, tied roots, Jordan blocks, and dependent age-specific reproduction remains open.



\end{document}